\documentclass[10pt]{amsart}
\usepackage{epsfig}
\usepackage{amsmath, amssymb, euscript,  enumerate}
\usepackage{bbm}
\usepackage{color}
\usepackage{verbatim}

\newcommand{\supp}{\text{\rm supp}}

\newcommand{\ap}{\alpha}             
\newcommand{\bt}{\beta}
\newcommand{\gm}{\gamma}             \newcommand{\Gm}{\Gamma}
\newcommand{\dt}{\delta}             
\newcommand{\vep}{\varepsilon}
\newcommand{\zt}{\zeta}

\newcommand{\ld}{\lambda}            \newcommand{\Ld}{\Lambda}
\newcommand{\sm}{\sigma}             
\newcommand{\vp}{\varphi}
\newcommand{\om}{\omega}             \newcommand{\Om}{\Omega}
\newcommand{\vr}{\varrho}            \newcommand{\iy}{\infty}
            
\newcommand{\f}{\frac}             \newcommand{\el}{\ell}

\newcommand{\fC}{{\mathfrak C}}

\newcommand{\fF}{{\mathfrak F}}

\newcommand{\fL}{{\mathfrak L}}

\newcommand{\ff}{{\mathfrak f}}
\newcommand{\fg}{{\mathfrak g}}
\newcommand{\fh}{{\mathfrak h}}

\newcommand{\fu}{{\mathfrak u}}
\newcommand{\fv}{{\mathfrak v}}

\newcommand{\fy}{{\mathfrak y}}

\newcommand{\BN}{{\mathbb N}}

\newcommand{\BR}{{\mathbb R}}

\newcommand{\cA}{{\mathcal A}}

\newcommand{\cC}{{\mathcal C}}

\newcommand{\cF}{{\mathcal F}}

\newcommand{\cI}{{\mathcal I}}
\newcommand{\cJ}{{\mathcal J}}
\newcommand{\cK}{{\mathcal K}}

\newcommand{\cM}{{\mathcal M}}

\newcommand{\cO}{{\mathcal O}}
\newcommand{\cP}{{\mathcal P}}
\newcommand{\cQ}{{\mathcal Q}}
\newcommand{\cR}{{\mathcal R}}
\newcommand{\cS}{{\mathcal S}}
\newcommand{\cT}{{\mathcal T}}

          \newcommand{\util}{{\tilde u}}

\newcommand{\la}{\langle}          \newcommand{\ra}{\rangle}
\newcommand{\s}{\setminus}         \newcommand{\ep}{\epsilon}
\newcommand{\n}{\nabla}            \newcommand{\e}{\eta}
\newcommand{\pa}{\partial}        
       
    \newcommand{\ds}{\displaystyle}

 \newcommand{\pf }{\noindent{\it Proof. }}
\newcommand{\rk }{\noindent{\bf Remark. }}
\newcommand{\bu }{\bar u} \newcommand{\bv }{\bar v}

\newcommand{\aee }{\text{\rm a.e.}} 
  \newcommand{\pv }{\text{\rm p.v.}}
  \newcommand{\dd }{\text{\rm d}}

\newcommand{\rL }{{\text{\rm L}}}  
   \newcommand{\rI}{{\text{\rm I}}}
  
\newcommand{\rX}{{\text{\rm X}}}  \newcommand{\rY}{{\text{\rm Y}}}

\newcommand{\esup}{{\text{\rm ess sup}}}
\newcommand{\spa}{{\text{\rm span}}}

\newtheorem{thm}[subsection]{Theorem}
\newtheorem{lemma}[subsection]{Lemma}
\newtheorem{cor}[subsection]{Corollary}

\newtheorem{defn}[subsection]{Definition}

\numberwithin{equation}{section}

\title[ Nonlocal Harnack inequalities ]{ Nonlocal Harnack inequalities \\ for nonlocal heat equations }
\author{ Yong-Cheol Kim }
\begin{document}
\begin{abstract} By applying the De Giorgi-Nash-Moser theory, we obtain nonlocal Harnack inequalities for weak solutions
of nonlocal parabolic equations given
by an integro-differential operator $\rL_K$ as follows;
\begin{equation*}\begin{cases} \rL_K u+\pa_t u=0 &\text{ in $\Om\times(-T,0]$ }\\
                u=g &\text{ in $\bigl((\BR^n\s\Om)\times (-T,0]\bigr)\cup\bigl(\Om\times\{t=-T\}\bigr)$ } \\
                              \end{cases}\end{equation*}
where $g\in C(\BR^n\times [-T,0])\cap L^{\iy}(\BR^n\times(-T,0])\cap H^s_T(\BR^n)$ and $\,\Om\,$ is a bounded domain in $\BR^n$
with Lipschitz
boundary. Moreover, we get nonlocal parabolic weak Harnack inequalities of the weak
solutions.
\end{abstract}

\address{$\bullet$ Yong-Cheol Kim : Department of Mathematics Education, Korea University, Seoul 136-701,
Republic of Korea }

\email{ychkim@korea.ac.kr}

\maketitle

\tableofcontents

\section{Introduction}

The study of fractional and nonlocal equations has recently been
done not only in pure mathematical analysis area but also in the
research that needs its concrete applications. The aim of this paper
is to obtain nonlocal Harnack inequalities for weak solutions of
nonlocal heat equations.

Let $\cK_0$ be the collection of all positive symmetric kernels
satisfying the uniformly ellipticity assumption
\begin{equation}\frac{(1-s)\lambda}{|y|^{n+2s}}\leq K(y)\leq
\frac{(1-s)\Lambda}{|y|^{n+2s}},\,\,0<s<1.
\end{equation}
Here the symmetricity means that $K(y)=K(-y)$ for all $y\in\BR^n$.
Then we consider the corresponding nonlocal operator $\rL_K$ given
by
\begin{equation}\rL_K u(x,t)=\pv\int_{\BR^n}\mu_t(u,x,y)K(y)\,dy
\end{equation}
where $\mu_t(u,x,y)=2 u(x,t)-u(x+y,t)-u(x-y,t)$. Set $\fL_0=\{\rL_K:
K\in\cK_0\}$. In particular, if $K(y)=c_{n,s}|y|^{-n-2s},\,s\in(0,1),$ where
$c_{n,s}$ is the normalization constant comparable to $s(1-s)$
given by
$$c_{n,s}=\f{1}{2}\int_{\BR^n}\f{1-\cos(\xi_1)}{|\xi|^{n+2s}}\,d\xi,$$
then $\rL_K=(-\Delta)^s$ is the fractional Laplacian and it
is well-known that $$\lim_{s\to 1^-}(-\Delta)^s u=-\Delta u$$
for any function $u$ in the Schwartz space $\cS(\BR^n)$.

In this paper, we study the boundary value problem for the following
nonlocal parabolic equations ${\bf NP}_{\Om_I}(f,g,h)$
\begin{equation}\begin{cases} \rL_K u+\pa_t u=f &\text{ in $\Om_I:=\Om\times I$ }\\
                u=g &\text{ in $(\BR^n\s\Om)\times I$ } \\
                u(x,-T)=h(x) &\text{ for $x\in\BR^n,$}
                                \end{cases}\end{equation}
where $I:=(-T,0]$ and $\Om$ is a bounded domain in $\BR^n$ with
Lipschitz boundary. More precisely speaking, by employing the De
Giorgi-Nash--Moser theory, we obtain nonlocal parabolic Harnack
inequalities for weak solutions of the nonlocal parabolic equation
${\bf NP}_{\Om_I}(0,g,g)$ where  $g\in C(\BR^n_{I_*})\cap
L^{\iy}(\BR^n_I)$ for $I_*:=[-T,0]$, and also we get nonlocal
parabolic weak Harnack inequalities of the weak solutions.

\,\,{\bf Notations.} We write the notations briefly for the readers as follows.

\,\noindent$\bullet$ For $r>0$ and $s\in(0,1)$, let us denote by $Q^0_r:=Q_r(x_0,t_0)=B_r^0\times I_{r,s}(t_0)$ and $Q_r=Q_r(0,0)$, where $B_r^0=B_r(x_0)$, $B_r=B_r(0)$ and $I_{r,s}(t_0)=(t_0-r^{2s},t_0]$.
Also, we denote by $\cQ_r^0=\cQ_r(x_0,t_0)=B_r^0\times(t_0-(2-\sm)r^{2s},t_0]$, $\cQ_r=\cQ_r(0,0)$,
\begin{equation*}\begin{split}
(Q^0_r)^+&:=Q^+_r(x_0,t_0)=B^0_r\times I_{r,s}^+(t_0),\,\,\, Q^+_r(0,0):=Q^+_r,\\
(Q^0_r)^-&:=Q^-_r(x_0,t_0)=B^0_r\times I_{r,s}^-(t_0),\,\,\,Q^-_r(0,0):=Q^-_r,
\end{split}\end{equation*}
where $\sm\in(0,1/2)$ is a constant to be given in Section 6,
$$I_{r,s}^+(t_0)=(t_0-\sm r^{2s},t_0]\text{ and }I_{r,s}^-(t_0)=\biggl(t_0-\biggl(\f{1}{2}+\sm\biggr)r^{2s},t_0-\f{1}{2}\,r^{2s}\biggr].$$
For simplicity, we denote by $I_{r,s}(0)=I_{r,s}$, $I^+_{r,s}(0)=I^+_{r,s}$ and  $I^-_{r,s}(0)=I^-_{r,s}$.

\,\noindent$\bullet$ For two quantities $a$ and $b$, we write $a\lesssim b$ (resp.
$a\gtrsim b$) if there is a universal constant $C>0$ ({\it depending
only on $\ld,\Ld,n,s$ and $\ep$}) such that $a\le C\,b$ (resp. $b\le
C\,a$).

\,\noindent$\bullet$ For $a,b\in\BR$, we denote by $a\vee b=\max\{a,b\}$ and
$a\wedge b=\min\{a,b\}$.

\,\noindent$\bullet$ Let $\cF^n_T$ and $\cF^n$ be the families of all real-valued
Lebesgue measurable functions on $\BR^n\times(-T,0]$ and $\BR^n$,
respectively. For $u\in\cF^n_T$, we write $[\fu(t)](x):=u(x,t)$ and
$[\rL_K\fu(t)](x)=\rL_K u(x,t)$. Let $H^s_T(\BR^n)$ denote the function space consisting of all functions $u\in\cF^n_T$ such that $\fu(t)\in H^s(\BR^n)$ for all $t\in I_*$.

\,\noindent$\bullet$ For $(x_0,t_0)\in\Om_I$ and $r>0$, we now define the {\it nonlocal parabolic tail} of the function $u$ in $Q^0_r\subset\Om_I$ by
\begin{equation}\cT_r(u;(x_0,t_0))=\f{2s}{|S^{n-1}|}\,\, r^{2s}\sup_{t\in I_{r,s}(t_0)}\int_{\BR^n\s B_r(x_0)}\f{|u(y,t)|}{|y-x_0|^{n+2s}}\,dy.
\end{equation}

The first Harnack-type inequality for globally nonnegative weak
solutions of local heat equations given on $\Om\times I$ was
obtained independently by Hadamard \cite{H} and Pini \cite{P}. After
this, the major influential contribution to the study in this
direction (in fact, for local parabolic equations of divergence type
given in $\Om\times I$) was made by J. Moser \cite{M1, M2}.
Interestingly, the phenomenon that the classical Harnack inequality
no longer works for nonlocal elliptic operators was recently
observed by Kassmann \cite{K1}. This unexpected fact motivated lots
of mathematicians to study the so-called {\it nonlocal Harnack
inequality}.

We now state our main results which are called {\it nonlocal Harnack
inequalities} and {\it weak Harnack inequalities} for weak solutions
of nonlocal heat equations as follows. Their proofs can be obtained from Theorem 7.4, 7.5, Corollary 7.6 and Appendix.

\begin{thm} Let $g\in C(\BR^n_{I_*})\cap L^{\iy}(\BR^n_I)\cap H^s_T(\BR^n)$.
If $\,u\in H^1(L;\rX_g(\Om))$ is a weak solution of the nonlocal parabolic equation ${\bf NP}_{\Om_I}(0,g,g)$ with $u\ge 0$ in $Q^0_R\subset\Om_I$, then there exists a constnat $c>0$ depending only on $n,s,\ld$ and $\Ld$ such that
\begin{equation}\sup_{(Q^0_r)^-}u\le c\,\inf_{(Q^0_r)^+}u+c\,\biggl(\f{r}{R}\biggr)^{2s}\cT_r(u^-;(x_0,t_0))
\end{equation} for any $r\in (0,R/5)$.
\end{thm}

We can easily obtain the following nonlocal parabolic Harnack inequalities for a nonnegative weak solutions of the nonlocal parabolic equation ${\bf NP}_{\Om_I}(0,g,g)$ as a natural by-product of Theorem 1.1. By the way, it is interesting that the result has no nonlocal parabolic tail. That means that the result coincides with that of local parabolic case.

\begin{cor} Let $g\in C(\BR^n_{I_*})\cap L^{\iy}(\BR^n_I)\cap H^s_T(\BR^n)$.
If $\,u\in H^1(L;\rX_g(\Om))$ is any nonnegative weak solution of the nonlocal parabolic equation ${\bf NP}_{\Om_I}(0,g,g)$, then there exists a constnat $c>0$ depending only on $n,s,\ld$ and $\Ld$ such that
\begin{equation}\sup_{(Q^0_r)^-}u\le c\,\inf_{(Q^0_r)^+}u
\end{equation} for any $r\in (0,R/5)$.
\end{cor}

\begin{thm} Let $g\in C(\BR^n_{I_*})\cap L^{\iy}(\BR^n_I)\cap H^s_T(\BR^n)$.
If $u\in H^1(L;\rX_g(\Om))$ is a weak solution of the nonlocal parabolic equation ${\bf NP}_{\Om_I}(0,g,g)$ with $u\ge 0$ in $Q^0_R\subset\Om_I$, then we have the estimate
\begin{equation}\biggl(\f{1}{2|(Q^0_r)^+|}\int_{(Q^0_r)^+}u^p\,dx\,dt\biggr)^{\f{1}{p}}\le \inf_{(Q^0_r)^+}u+\f{4}{3}\,\biggl(\f{r}{R}\biggr)^{2s}\cT_r(u^-;(x_0,t_0))
\end{equation} for any $p\in(0,1)$ and $r\in (0,R)$.
\end{thm}

\rk (a) In case that $\Om=\BR^n$, using the De Giorgi method, Caffarelli, Chan and Vasseur proved in \cite{CCV} that any weak solution to the equation (1.3) with initial data $h$ in $L^2(\BR^n)$ is uniformly bounded and H$\ddot {\rm o}$lder continuous.

\,(b) The ellitic result of this problem was obtained by Di Castro, Kuusi and Palatucci \cite{DKP}. As a matter of fact, when $p\in(1,\iy)$, they proved nonlocal Harnack inequalities for elliptic nonlocal $p$-Laplacian equations there. Also, they obtained H$\ddot {\rm o}$lder regularity in \cite{DKP1}.

\,(c) Using the Moser's iteration method, Felsinger and Kassmann obtained weak parabolic Harnack inequality and  H$\ddot {\rm o}$lder regularity in \cite{FK}. Also, despite of failure for getting the classical Harnack inequality as mentioned above, the first attempt to obtain the nonlocal Harnack inequalities of the form (1.5) for the fractional elliptic equations was tried by Kassmann (see \cite{K2}.

\,(d) The H${\rm {\ddot o}}$lder continuity of weak solutions of the nonlocal parabolic equation ${\bf NP}_{\Om_I}(0,g,g)$ with $u\ge 0$ in $Q^0_R\subset\Om_I$ and $g\in C(\BR^n_{I_*})\cap L^{\iy}(\BR^n_I)\cap H^s_T(\BR^n)$ was obtained in \cite{K}.

\,(e) In \cite{BBK}, Barlow, Bass and Kumagai gave a probabilistic proof for parabolic Harnack inequality by using the connection between stochastic processes and equations similar to the nonlocal equation (1.3).

\,(f) In \cite{BSV}, Bonforte, Sire and V\'azquez established an optimal existence and uniqueness theory for the Cauchy problem for the fractional heat equations given in $\BR^n\times I$.

\,\, The paper is organized as follows. In Section 2, we furnish the
function spaces and the definition of weak solutions of the nonlocal
parabolic equations given in (1.3), and also give two well-known
useful lemmas. The maximum principle and comparison principle of
weak solutions for the nonlocal heat equations are obtained in
Section 3. In Section 4, we obtain a relation between weak solutions
and viscosity solutions of the nonlocal heat equations, which makes
its weak solutions possible to enjoy the previous nice results on
its viscosity solutions. In Section 5, we get nonlocal weak Harnack
inequality for the nonlocal heat equation which is useful in proving
its nonlocal parabolic Harnack inequality. It turned out that, in
the elliptic case, Poincar\'e inequality was one of the crucial
tools for the proof of classical Harnack inequality and no longer
depends on the given partial differential equations. However, the
fractional Poincar\'e inequality in the parabolic sense is not
available for a general weak solution $u\in L^2(I,\rX_0(\Om))$. In
Section 6, we obtain parabolic fractional Poincar\'e inequality
depending on the nonlocal heat equations. In Section 7, we establish
the proof of nonlocal parabolic Harnack inequality by applying the De Giorgi-Nash-Moser theory \cite{D,N,M}. Finally, in Appendix, we give the proof of the existence and uniqueness for weak
solutions of the nonlocal heat equations which is based on the
results for the weak formulation of the nonlocal eigenvalue problem
of elliptic type \cite{SV} and the Galerkin's method.

\section{Preliminaries}
Let $\rY$ be a real Banach space with norm $\|\cdot\|$
and let $\cF^\rY_T$ be the family of all measurable vector-valued
functions $\fu:I\to\rY$. For $1\le p\le\iy$, we introduce
vector-valued function spaces
$L^p(I;\rY)=\{\fu\in\cF^\rY_T:\|\fu\|_{L^p(I;\rY)}<\iy\},$
where
\begin{equation}\begin{split}\|u\|_{L^p(I;\rY)}&:=\biggl(\int^0_{-T}\|\fu(t)\|^p\,dt\biggr)^{1/p}
\,\,\,\text{ for $1\le p<\iy$, } \\
\|u\|_{L^{\iy}(I;\rY)}&:=\esup_{t\in I}\|\fu(t)\|\,\,\,\text{ for
$p=\iy$. }
\end{split}\end{equation}
We also consider the function space $C(I;\rY)$ consisting of all functions $u\in\cF^n_T$ such that $\fu:I\to\rY$ a continuous vector-valued function satisfying
\begin{equation}\|u\|_{C(I;\rY)}:=\sup_{t\in I}\|\fu(t)\|<\iy.
\end{equation}

Let $\fu\in L^1(I;\rY)$ Then we say that $\fv\in L^1(I;\rY)$ is the
{\it weak derivative} of $\fu$ and we write $\fv=\fu'$ if
$$\int^0_{-T}\vp'(t)\fu(t)\,dt=-\int^0_{-T}\vp(t)\fv(t)\,dt$$ for all
testing functions $\vp\in C^{\iy}_c(I)$.

Let $\Om$ be a bounded domain in $\BR^n$ with Lipschitz boundary
and let $K\in\cK_0$. Let $\rX(\Om)$ be the linear function space of all
Lebesgue measurable functions $v\in\cF^n$ such that $v|_\Om\in L^2(\Om)$
and
\begin{equation*}
\ds\iint_{\BR^{2n}_\Om}\f{|v(x)-v(y)|^2}{|x-y|^{n+2s}}\,dx\,dy<\iy
\end{equation*}
where $\BR^{2n}_\Om:=\BR^{2n}\s(\Om^c\times \Om^c)$. We
also set
\begin{equation}\rX_0(\Om)=\{v\in\rX(\Om):v=0\,\,\aee\text{ in $\BR^n\s\Om$ }\}
\end{equation}
Since $C^2_0(\Om)\subset\rX_0(\Om)$, we see that $\rX(\Om)$ and $\rX_0(\Om)$ are not
empty. Then we see that $(\rX(\Om),\|\cdot\|_{\rX(\Om)})$ is a normed space,
where the norm $\|\cdot\|_{\rX(\Om)}$ defined by
\begin{equation}\|v\|_{\rX(\Om)}:=\|v\|_{L^2(\Om)}+\biggl(\iint_{\BR^{2n}_\Om}\f{|v(x)-v(y)|^2}{|x-y|^{n+2s}}\,dx\,dy\biggr)^{1/2}<\iy,
\,\,\,v\in\rX(\Om).
\end{equation}
We denote by $H^s(\Om)$ the usual fractional Sobolev space
with the norm
\begin{equation}\|v\|_{H^s(\Om)}:=\|v\|_{L^2(\Om)}+[v]_{H^s(\Om)}<\iy
\end{equation} where the seminorm $[\,\cdot\,]_{H^s(\Om)}$ is defined by
$$[v]_{H^s(\Om)}:=[v]_{W^{s,2}(\Om)}=\biggl(\iint_{\Om\times\Om}
\f{|v(x)-v(y)|^2}{|x-y|^{n+2s}}\,dx\,dy\biggr)^{1/2}.$$
Then it is well-known \cite{SV} that
\begin{equation}\|v\|_{H^s(\Om)}\le\|v\|_{H^s(\BR^n)}\le
c(\ld,s)\|v\|_{\rX(\Om)}
\end{equation}
for any $v\in\rX_0(\Om)$, where
$c(\ld,s)=\max\{1,[\ld(1-s)]^{-1/2}\}$. Also, there is a
constant $c_0>1$ depending only on $n,\ld,\sm$ and $\Om$ such that
\begin{equation}\|v\|^2_{\rX_0(\Om)}\le\|v\|^2_{\rX}\le
c_0\|v\|^2_{\rX_0(\Om)}
\end{equation}
for any $v\in\rX_0(\Om)$; that is,
\begin{equation}\|v\|_{\rX_0(\Om)}:=\biggl(\iint_{\BR^{2n}_\Om}\f{|v(x)-v(y)|^2}{|x-y|^{n+2s}}\,dx\,dy\biggr)^{1/2}
\end{equation} is a norm on $\rX_0(\Om)$ equivalent to (2.4). Moreover it is known that $(\rX_0(\Om),\|\cdot\|_{\rX_0(\Om)})$ is a
Hilbert space with inner product
\begin{equation}\la
u,v\ra_{\rX_0(\Om)}:=\iint_{\BR^{2n}_\Om}(u(x)-u(y))(v(x)-v(y)\,d_K(x,y)
\end{equation} where $d_K(x,y):=K(x-y)\,dx\,dy$.
Let $\rX_0^*(\Om)$ be the dual space of $\rX_0(\Om)$, i.e. the family of all
bounded linear functionals on $\rX_0(\Om)$. Then we see that $\rX_0(\Om)\subset\rX(\Om)\subset L^2(\Om)\subset\rX_0^*(\Om)$ and
$(\rX_0^*(\Om),\|\cdot\|_{\rX_0^*(\Om)})$ is a normed space, where the norm
$\|\cdot\|_{\rX_0^*(\Om)}$ is given by
$$\|u\|_{\rX_0^*(\Om)}:=\sup\{u(v):v\in\rX_0(\Om), \|v\|_{\rX_0(\Om)}\le
1\},\,\,u\in\rX_0^*(\Om).$$

In what follows, for a Banach space $(B,\|\cdot\|_B)$ with its dual space $(B^*,\|\cdot\|_{B^*})$, we consider a vector-valued Banach space
\begin{equation*}
H^1(I;B)=\{u\in L^2(I;B):\fu'\in L^2(I;B^*)\}
\end{equation*}
with the norm
\begin{equation*}\|u\|_{H^1(I;B)}=\biggl(\int^0_{-T}\|\fu(t)\|^2_B\,dt\biggr)^{\f{1}{2}}+\biggl(\int^0_{-T}\|\fu'(t)\|^2_{B^*}\,dt\biggr)^{\f{1}{2}}<\iy.
\end{equation*}
For $g\in H^s(\BR^n)$, we consider the convex subsets of $H^s(\BR^n)$ by \begin{equation*}\begin{split}\rX_g^\pm(\Om)&=\{v\in H^s(\BR^n):(g-v)^{\pm}\in\rX_0(\Om)\},\\
\rX_g(\Om)&:=\rX_g^+(\Om)\cap\rX_g^-(\Om)=\{v\in H^s(\BR^n):g-v\in\rX_0(\Om)\}.
\end{split}\end{equation*}
For $g\in H^s_T(\BR^n)$, we define the {\it convex subsets} of the space $H^s_T(\BR^n)$ by
\begin{equation*}\begin{split}H^1(I;\rX^{\pm}_g(\Om))&=\{u\in H^s_T(\BR^n):(g-u)_{\pm}\in H^1(I;\rX_0(\Om))\},\\
H^1(I;\rX_g(\Om))&:=H^1(I;\rX^+_g(\Om))\cap H^1(I;\rX^-_g(\Om))\\
&=\{u\in H^s_T(\BR^n):u-g\in H^1(I;\rX_0(\Om))\}.
\end{split}\end{equation*}

\,\rk If $\fu\in H^1(I;\rX_0(\Om))$, then it is well-known that
(a) $\fu\in C(I;L^2(\Om))$ after being modified on a set of measure
zero, (b) the function $\ap$ defined by $\ap(t)=\|\fu(t)\|^2_{L^2(\Om)}$
is absolutely continuous, and moreover
$\ap'(t)=2\la\fu'(t),\fu(t)\ra_{L^2(\Om)}$ for $\aee$ $t\in I$, and (c) there is a constant $C>0$ depending only on $T$ such that
$$\sup_{t\in I}\|\fu(t)\|_{L^2(\Om)}\le
C\,\|\fu\|_{H^1(I;\rX_0(\Om))}.$$

\,\,In order to define weak solutions, we consider a bilinear form defined by
$$\la u,v\ra_K=\iint_{\BR^n\times\BR^n}(u(x)-u(y))(v(x)-v(y))\,K(x-y)\,dx\,dy\,\,\text{ for $u,v\in\rX(\Om)$.}$$

\begin{defn} Let $g\in H^s_T(\BR^n)$ and $f\in L^2(I;\rX^*_0(\Om))$. Then we say that a function $u\in H^1(I;\rX^-_g(\Om))$ $(\,u\in H^1(I;\rX^+_g(\Om))\,)$ is a {\rm weak subsolution (\,weak supersolution\,)} of the nonlocal
parabolic equation ${\bf NP}_{\Om_I}(f,g,h)$ given in $(1.3)$, if it satisfies
\begin{equation}\la\fu(t),\vp\ra_K+\la\fu'(t)-\ff(t),\vp\ra\le 0\,\,(\,\ge 0\,)\end{equation}
for any nonnegative $\vp\in\rX_0(\Om)$ and $\aee$ $t\in I$, and
\begin{equation}\fu(-T)=h,\end{equation}
where $\la\cdot,\cdot\ra$ denotes the dual pair between $\rX_0(\Om)$ and $\rX^*_0(\Om)$.
Also, we say that a function $u$ is a {\rm weak solution} of the equation ${\bf NP}_{\Om_I}(f,g,h)$, if it is both a weak subsolution and a weak supersolution, i.e.
\begin{equation}\la\fu(t),\vp\ra_K+\la\fu'(t)-\ff(t),\vp\ra=0\end{equation}
for any $\vp\in\rX_0(\Om)$ and $\aee$ $t\in I$, and $\fu(-T)=h.$
\end{defn}

\,\,In order to prove our results, we need two well-known lemmas to be given in the following (see \cite{GT}).

\begin{lemma} Let $\{N_k\}_{k=0}^{\iy}\subset\BR$ be a sequence of positive numbers such that
$$N_{k+1}\le d_0\,e_0^k N_k^{1+\e}$$
where $d_0,\e>0$ and $e_0>1$. If $N_0\le d_0^{-1/\e} e_0^{-1/\e^2}$, then we have that $N_k\le e_0^{-k/\e}\,N_0$ for any $k=0,1,\cdots$ and moreover $\lim_{k\to\iy}N_k=0$.
\end{lemma}

\begin{lemma} Let $f$ be a nonnegative bounded function defined in $[t_0,t_1]$, where $0\le t_0<t_1$. Suppose that there are nonnegative constants $c_1,c_2,\theta,$ and $\e\in(0,1)$ such that
$$f(t)\le c_1(\tau-t)^{-\theta}+c_2+\e\,f(\tau)$$ for any $t,\tau\in[t_1,t_2]$ with $t<\tau$. Then there exists a constant $c>0$ $($depending only on $\theta$ and $\e$\,$)$ such that
$$f(\rho)\le c[c_1(R-\rho)^{-\theta}+c_2]$$ for any $\rho,R\in[t_1,t_2]$ with $\rho<R$.
\end{lemma}

\section{Maximum principle and comparison principle}

In this section, we furnish Maximum Principle and Comparison Principle for weak solutions of the nonlocal parabolic equations ${\bf NP}_{\\\cO_J}(0,g,g)$ where $\cO\subset\BR^n$ is a bounded open set and $J:=[a,b)\subset I$ is a half-open interval. We denote by $\BR^n_{I_*}:=\BR^n\times I_*$ for $I_*=[-T,0]$.

\begin{lemma}  If $u\in H^1(J;\rX^-_g(\cO))$ is a weak subsolution of the nonlocal parabolic equation ${\bf NP}_{\cO_J}(0,g,g)$ given in $(1.3)$ and $u=g\le 0$ in  $\BR^n_{I_*}\s\cO_J$, then $u\le 0$ in $\BR^n\times I$.
\end{lemma}

\pf By the assumption, we see that $u^+=0$ in $(\BR^n\s\cO)\times J$, and thus $u^+\in H^1(J;\rX_0(\cO))$. Thus we can use $u^+$ as  a testing function in the weak formulation.
Observing that $u^+(x,t)u^-(x,t)=0$ and $u^+(x,t)u^-(y,t)\ge 0$ for $\aee\,x,y\in\BR^n$ and $t\in J$, it follows from the fractional sobolev inequality that
\begin{equation*}\begin{split}
0&\ge\la\fu(t),\fu^+(t)\ra_{\rX_0(\cO)}+\la\pa_t\fu(t),\fu^+(t)\ra_{L^2(\cO)}\\
&=\|\fu^+(t)\|^2_{\rX_0(\cO)}-\la\fu^-(t),\fu^+(t)\ra_{\rX_0(\cO)}\\
&\qquad\qquad\qquad\qquad+\pa_t\bigl(\la \fu^+(t),\fu^+(t)\ra_{L^2(\cO)}-\la\fu^-(t),\fu^+(t)\ra_{L^2(\cO)}\bigr)\\
&\ge c\|\fu^+(t)\|^2_{L^2(\cO)}+\pa_t\|\fu^+(t)\|^2_{L^2(\cO)}\\
&\qquad\qquad\qquad\qquad-\iint_{\BR^{2n}_{\cO}}\f{(u^-(x,t)-u^-(y,t))(u^+(x,t)-u^+(y,t))}{|x-y|^{n+2s}}\,dx\,dy\\
&\ge c\|\fu^+(t)\|^2_{L^2(\cO)}+\pa_t\|\fu^+(t)\|^2_{L^2(\cO)}\\
&\qquad\qquad\qquad\qquad+\iint_{\BR^{2n}_{\cO}}\f{u^-(x,t)u^+(y,t)+u^+(x,t)u^-(y,t)}{|x-y|^{n+2s}}\,dx\,dy\\
&\ge c\|\fu^+(t)\|^2_{L^2(\cO)}+\pa_t\|\fu^+(t)\|^2_{L^2(\cO)}.
\end{split}\end{equation*}
Thus, by Gronwall's inequality and the assumption on the initial values, we conclude that
\begin{equation*}\|\fu^+(t)\|^2_{L^2(\cO)}\le e^{-c(t+T)}\|\fu^+(-T)\|^2_{L^2(\cO)}=0.
\end{equation*}
This implies that $u^+=0$ in $\cO\times J$. Hence we are done. \qed

\begin{cor}  If $u\in H^1(J;\rX^+_g(\cO))$ is a weak supersolution of the nonlocal parabolic equation ${\bf NP}_{\cO_J}(0,g,g)$ given in $(1.3)$ and $u=g\ge 0$ in  $\BR^n_{I_*}\s\cO_J$, then $u\ge 0$ in $\BR^n\times I$.
\end{cor}

\begin{cor} If $u\in H^1(I;\rX^-_{g_1}(\cO))$ is a weak subsolution of the nonlocal parabolic equation ${\bf NP}_{\cO_J}(0,g_1,g_1)$ and $v\in H^1(I;\rX^+_{g_2}(\cO))$ is a weak supersolution of the nonlocal parabolic equation ${\bf NP}_{\cO_J}(0,g_2,g_2)$ such that $g_1\le g_2$ in  $\BR^n_{I_*}\s\cO_J$, then $u\le v$ in $\BR^n\times I$.
\end{cor}

\section{Weak and viscosity solutions }

In this section, we get boundedness and continuity on $\BR^n_I$ of weak solutions of the nonlocal parabolic equation ${\bf NP}_{\Om_I}(0,g,g)$ with boundary condition $g\in C(\BR^n_{I_*})\cap L^{\iy}(\BR^n_I)$ where $\BR^n_{I_*}:=\BR^n\times I_*$ for $I_*=[-T,0]$ and we study a relation between weak solutions and viscosity solutions of the equation ${\bf NP}_{\Om_I}(0,g,g)$. The latter one makes its weak solutions possible to enjoy the previous nice results
on its viscosity solutions.

Let us define viscosity solutions. Let $\cP(\BR^{n+1})$ denote the class of all parabolic quadratic polynomials of the form
$$p(x,t)=\sum_{i=1}^n\sum_{j=1}^n a_{ij}x_i x_j+\sum_{i=1}^n b_i x_i +ct +d,$$ where $a_{ij},b_i,c,d\in\BR.$ A upper (lower) semicontinuous function $u:\BR^n\times I\to\BR$ is called a {\it viscosity subsolution} ({\it viscosity supersolution}) of the equation ${\bf NP}_{\Om_I}(0,g,g)$ and we write $\rL_K u+\pa_t u\le 0$ (res. $\rL_K u+\pa_t u\ge 0$) on $\Om_I$ in the viscosity sense, if for each $(x,t)\in\Om_I$ there is a neighborhood $Q_r(x,t)\subset\Om_I$ such that  $\rL_K u+\pa_t u$ is well-defined and  $\rL_K v(x,t)+\pa_t p(x,t)\le 0$ (res. $\rL_K v(x,t)+\pa_t p(x,t)\ge 0$) for $v=p\mathbbm{1}_{Q_r(x,t)}+u\mathbbm{1}_{\BR^n_I\s Q_r(x,t)}$ whenever $p\in\cP(\BR^{n+1})$ with $p(x,t)=u(x,t)$ and $p>u$ (res. $p<u$) on $Q_r(x,t)\s\{(x,t)\}$ exists.
Moreover, a function $u:\BR^n\times I\to\BR$ is called a {\it viscosity solution}, if it is both a viscosity subsolution and a viscosity supersolution of the equation.

\begin{thm} If $u\in H^1(I;\rX_g(\Om))$ is a weak solution of the nonlocal parabolic equation ${\bf NP}_{\Om_I}(0,g,g)$ for $g\in C(\BR^n_{I_*})\cap L^{\iy}(\BR^n_I)$,
then $u\in L^{\iy}(\BR^n_I)$.
\end{thm}

\pf By multiplying $u$ by a sufficiently small constant, we may assume that
\begin{equation}\|u\|^2_{L^2(\Om_I)}\le\dt
\end{equation} where $\dt>0$ is some small constant to be determined later. Let $M=2\|g\|_{L^{\iy}(\BR^n_I)}$ and take any $a,b\in(0,1)$ so that $a+b=1$.
For $k\in\BN$, let $M_k=M(1-2^{-k})$, $u_k=u-M_k$, $w_k=u_k^+$ and $N_k=\|w_k\|^2_{L^2(\Om_{I_k})}$, where $I_k=(-T_k,0]$ for  $T_k=(a+2^{-k}b)T$ and $I_0=I$ for $T_0=T$.
Then we see that, for any $k\in\BN$, $M_{k+1}>M_k$ and $u_{k+1}<u_k$ , and so $w_{k+1}\le w_k$. Moreover, on $(\BR^n\s\Om)\times I$, we have that
$$u_{k+1}=g-M+2^{-k-1}M\le-\f{M}{2}+2^{-k-1}M\le 0\,\text{ for all $k=0,1,\cdots$\,.}$$
So we have that $w_{k+1}=0$ on $(\BR^n\s\Om)\times I$. We now use $\vp_{k+1}=w_{k+1}\e_{k+1}^2$ as a testing function in the weak formulation of the equation, where $\e_{k+1}\in C_c^{\iy}(-T_k,\iy)$ is a function such that $0\le\e_{k+1}\le 1$ in $\BR$, $\e_{k+1}=1$ in $(-T_{k+1},\iy)$ and $0\le\e'_{k+1}\le 2^{k+2}(bT)^{-1}$ in $\BR$. Then we have that
\begin{equation*}\begin{split}
\int_{-T_k}^\tau\int_{\Om}(\pa_t u)\vp_{k+1}\,dx\,dt+\rI(u,\vp_{k+1})=0
\end{split}\end{equation*} for any $\tau\in(-T_k,0]$,
where the bilinear operator is given by
\begin{equation*}\begin{split}&\rI(u,\vp_{k+1})\\&=\int_{-T_k}^\tau\e_{k+1}^2(t)
\iint_{\Om\times\Om}(u(x,t)-u(y,t))(\vp_{k+1}(x,t)-\vp_{k+1}(y,t))\,d_K(x,y)\,dt.
\end{split}\end{equation*}
The first term in the left-hand side of the above equality can be evaluated by
\begin{equation}\begin{split}
&\int_{-T_k}^\tau\int_{\Om}(\pa_t u)\vp_{k+1}\,dx\,dt=\int_{-T_k}^\tau
\int_{\Om}\bigl[\pa_t(w_{k+1}^2\e_{k+1}^2)-2 w_{k+1}^2\e_{k+1}\e'_{k+1}\bigr]dx\,dt\\
&\qquad=\e_{k+1}^2(\tau)\int_{\Om}w_{k+1}^2(x,\tau)\,dx-2\int_{-T_k}^\tau\e_{k+1}(t)\e'_{k+1}(t)\int_{\Om}w_{k+1}^2(x,t)\,dx\,dt.
\end{split}\end{equation}
We next split $\rI(u,\vp_{k+1})$ into two parts as follows;
\begin{equation}\begin{split}&\rI(u,\vp_{k+1})\\&=\int_{-T_k}^\tau\e_{k+1}^2(t)\iint_{\Om\times\Om}(u(x,t)-u(y,t))(\vp_{k+1}(x,t)-\vp_{k+1}(y,t))d_K(x,y)dt\\
&\qquad\quad+2\int_{-T_k}^\tau\e_{k+1}^2(t)\int_{\BR^n\s\Om}\int_{\Om}(u(x,t)-u(y,t))\,\vp_{k+1}(x,t)\,d_K(x,y)\,dt\\
&:=I_1+2 I_2.
\end{split}\end{equation}
For the estimate of $I_1$, we first observe that
\begin{equation}\begin{split}&(u(x,t)-u(y,t))(\vp_{k+1}(x,t)-\vp_{k+1}(y,t))\\
&\qquad\qquad\qquad\ge(w_{k+1}(x,t)-w_{k+1}(y,t))(\vp_{k+1}(x,t)-\vp_{k+1}(y,t))
\end{split}\end{equation}
whenever $(x,t),(y,t)\in\Om_I$; indeed, without loss of generality we may assume that $u(x,t)\ge u(y,t)$. Then it can easily be checked by considering two possible cases (i) $u(x,t)\ge u(y,t)>M_{k+1}$,  (ii) $u(x,t)>M_{k+1}$, $u(y,t)\le M_{k+1}$, and (iii) $M_{k+1}\ge u(x,t)\ge u(y,t)$.
For the estimate of $I_2$, we note that
\begin{equation*}\begin{split}
(u(x,t)-u(y,t))\vp_{k+1}(x,t)&\ge-(u(y,t)-u(x,t))_+(u(x,t)-M_{k+1})_+\e_{k+1}^2(t)\\
&\ge-(u(y,t)-M_{k+1})_+(u(x,t)-M_{k+1})_+\e_{k+1}^2(t)\\
&=-w_{k+1}(y,t)w_{k+1}(x,t)\e_{k+1}^2(t)
\end{split}\end{equation*}
and thus we have that
\begin{equation}\begin{split}
I_2&\ge-\int_{-T_k}^\tau\e_{k+1}^2(t)\int_{\BR^n\s\Om}\int_{\Om}w_{k+1}(y,t)w_{k+1}(x,t)\,d_K(x,y)\,dt=0.
\end{split}\end{equation}
Since the following equality
\begin{equation}\begin{split}&(w_{k+1}(x,t)-w_{k+1}(y,t))(\vp_{k+1}(x,t)-\vp_{k+1}(y,t))\\&\qquad\qquad\qquad\qquad=(w_{k+1}(x,t)-w_{k+1}(y,t))^2\,\e_{k+1}^2(t)\\
\end{split}\end{equation}
is always true, it follows from (4.3), (4.4), (4.5) and (4.6) that
\begin{equation}\begin{split}
&\rI(u,\vp_{k+1})\ge\int_{-T}^\tau\e_{k+1}^2(t)\iint_{\Om\times\Om}(w_{k+1}(x,t)-w_{k+1}(y,t))^2\,d_K(x,y)\,dt.
\end{split}\end{equation}
Thus it follows from (4.2) and (4.7) that
\begin{equation}\begin{split}
&\sup_{t\in(-T_k,0]}\e_{k+1}^2(t)\|w_{k+1}(\cdot,t)\|^2_{L^2(\Om)}\\
&\qquad+\int_{-T_k}^{0}\e_{k+1}^2(t)\iint_{\Om\times\Om}(w_{k+1}(x,t)-w_{k+1}(y,t))^2\,d_K(x,y)\,dt\\
&\quad\le 2\int_{-T_k}^{0}\e_{k+1}(t)\e'_{k+1}(t)\|w_{k+1}(\cdot,t)\|^2_{L^2(\Om)}\,dt.
\end{split}\end{equation}
Applying a well-known parabolic version of the fractional Sobolev inequality to (4.8), we obtain that
\begin{equation}\begin{split}
&\bigl\|w_{k+1}\bigr\|^{2\ap}_{L^{2\ap}(\Om_{I_{k+1}})}\le\iint_{\Om_{I_k}}|\e_{k+1}(t)w_{k+1}(x,t)|^{2\ap}\,dx\,dt\\
&\qquad\le c\,(2\|\e'_{k+1}\|_{L^{\iy}(\BR)}+1)^{\ap}\biggl(\int_{-T_k}^0\e_{k+1}(t)\int_{\Om}w_{k+1}^2(x,t)\,dx\,dt\biggr)^{\ap}\\
&\qquad\le c\,\f{2^{k\ap}}{(bT)^{\ap}}\,\bigl\|w_k\bigr\|^{2\ap}_{L^2(\Om_{I_k})},
\end{split}\end{equation} where $\ap=1+\f{2s}{n}$.
Since $\{w_{k+1}>0\}\subset\{w_k>2^{-(k+1)}M\}$, we have that
\begin{equation}\begin{split}
N_k&\ge\iint_{\Om_{I_k}\cap\{w_k>2^{-(k+1)}M\}}w_k^2(x,t)\,dx\,dt\\
&\ge 2^{-2(k+1)}M^2\,|\Om_{I_k}\cap\{w_k>2^{-(k+1)}M\}|\\
&\ge 2^{-2(k+1)}M^2\,|\Om_{I_k}\cap\{w_{k+1}>0\}|.\\
\end{split}\end{equation}
Thus, by (4.9), (4.10) and H$\rm\ddot o$lder's inequality, we have that
\begin{equation*}\begin{split}N_{k+1}&\le\biggl(\iint_{\Om_{I_{k+1}}}w^{2\ap}_{k+1}(x,t)\,dx\,dt\biggr)^{1/\ap}\,|\Om_{I_k}\cap\{w_k>0\}|^{\f{2s}{n+2s}}\\
&\le\f{c}{bT M^{\f{4s}{n+2s}}}\,2^{(1+\f{4s}{n+2s})k}\,N_k^{1+\f{2s}{n+2s}}.
\end{split}\end{equation*}
Set $\ds d_0=\f{c}{bT M^{\f{4s}{n+2s}}}>1$, $\ds e_0=2^{1+\f{4s}{n+2s}}>1$ and $\ds\e=\f{4s}{n+2s}>0$. If we select $\dt>0$ so that
$$\dt\le d_0^{-1/\e} e_0^{-1/\e^2},$$ then by (4.1) we have that $N_0=\|u^+\|^2_{L^2(\Om_I)}=\|u\|^2_{L^2(\Om_I)}\le d_0^{-1/\e} e_0^{-1/\e^2}$. Thus, by Lemma 2.2, we have that
\begin{equation*}0=\lim_{k\to\iy}N_k=\|(u-M)_+\|^2_{L^2(\Om_{(-aT,0]})},
\end{equation*} and hence $u\le M$ in $\Om_{(-aT,0]}$.
Also, by applying $-u$ instead of $u$, we have that $u\ge-M$ in $\Om_{(-aT,0]}$. Thus we have that
$|u|\le M$ in $\Om_{(-aT,0]}$ for any $a,b\in(0,1)$ with $a+b=1$.
Therefore, taking $a\downarrow 1$, we obtain that
$$\|u\|_{L^{\iy}(\Om_I)}\le M,$$
and we conclude that
$$\|u\|_{L^{\iy}(\BR^n_I)}\le 2\|g\|_{L^{\iy}(\BR^n_I)}<\iy.$$
Hence we complete the proof. \qed

\,\, In the next theorem, we obtain the global continuity of weak solutions of the nonlocal heat equations with certain boundary condition whose proof is based on the idea of that of the elliptic case \cite{SV1}.

\begin{thm} If $u\in H^1(I;\rX_g(\Om))$ is a weak solution of the nonlocal parabolic equation ${\bf NP}_{\Om_I}(0,g,g)$ for $g\in C(\BR^n_{I_*})\cap L^{\iy}(\BR^n_I)$,
then $u\in C(\BR^n_I)$.
\end{thm}

\pf For a contrapositive proof, we assume that there exists some $(x_0,t_0)\in\BR^n_{I_*}$ and sequence $(x_k,t_k)\in\BR^n_{I_*}$ with $\lim_{k\to\iy}(x_k,t_k)=(x_0,t_0)$ such that
$$|u(x_k,t_k)-u(x_0,t_0)|\ge\e_0$$ for some $\e_0>0$. Without loss of generality, we may assume that
\begin{equation}u(x_k,t_k)-u(x_0,t_0)\ge\e_0\,;
\end{equation}
for, the other case can be done in a similar way. Then we first claim that
\begin{equation} (x_0,t_0)\in\pa_p\Om_I.
\end{equation}
Indeed, $(x_0,t_0)\in\overline\Om\times I_*$, because $u$ is continuous in $\BR^n_{I_*}\s\Om_I$. Moreover, it is impossible that $(u_0,t_0)\in\Om\times I$, because we see from the local interior regularity results [FK] that $u$ is continuous in any compact domain contained in $\Om_I$. This implies (4.12).

For $\vep\in(0,1]$, let $\Om^{\vep}$ be a smooth $\vep$-neighborhood of $\Om$, i.e. a set with smooth boundary such that
\begin{equation}\bigcup_{x\in\Om}B_{\vep/2}(x)\subset\Om^{\vep}\subset\bigcup_{x\in\Om}B_{\vep}(x).
\end{equation}
If we consider the function $$g_\vep(x,t)=(1-\vp_\vep(x))u(x,t)+\vp_\vep(x)(t+2T)$$ where the function $\vp_\vep\in C^{\iy}(\BR^n)$ satisfies $\vp_\vep=1$ in $\Om_{\vep/4}$ and $\vp_\vep=0$ in $\BR^n\s\Om_{\vep/2}$, then we see that $g_\vep\in C(\BR^n_{I_*})$ and $g_\vep=u$ in $\BR^n_{I_*}\s\Om^{\vep}_I$.
Let $u_{\vep}\in H^1(I;\rX_{g_{\vep}}(\Om^{\vep}))$ be a weak solution of the nonlocal parabolic equation
\begin{equation}\begin{cases}\rL_K u_{\vep}+\pa_t u_{\vep}=0 & \text{ in $\Om^{\vep}_I$,}\\
u_{\vep}=g_\vep=u & \text{ in $\BR^n_{I_*}\s\Om^{\vep}_I$.}
\end{cases}\end{equation}
Then we see that $u_{\vep}\in L^{\iy}(\BR^n_I)$.
From Theorem 1.2 [FK], there are a constant $\ap\in(0,\sm_0)$ (depending only on $\ld,\Ld,n$ and $\sm_0\in(0,2)$) such that
\begin{equation}[u_{\vep}]_{C^{\ap}(Q)}\le\f{\|u_{\vep}\|_{L^{\iy}(\BR^n_I)}}{\e^{\ap}_{\vep}}
\end{equation}
for any $Q\Subset\Om^{\vep}_I$, where $\e_{\vep}=\e_{\vep}(Q)>0$ is some constant depending on $Q$ and $\Om^{\vep}_I$.
We consider the modulus of continuity $\vr$ of $u_{\vep}$ in $\Om^{\vep}_I$ defined by
\begin{equation*}\vr(\bt)=\sup_{Q\Subset\Om^{\vep}_I}\sup_{\substack{(x,t),(y,\tau)\in Q\\(|x-y|^{2s}+|t-\tau|)^{\ap/2s}<\e^{\ap}_{\vep}\bt/\|u_{\vep}\|_{L^{\iy}(\BR^n_I)}}}|u_{\vep}(x,t)-u_{\vep}(y,\tau)|.
\end{equation*}
Here we note that $\vr$ no longer depend on $\vep$; indeed, if $(x,t),(y,\tau)\in Q$ for some $Q\Subset\Om^{\vep}_I$ with $(|x-y|^{2s}+|t-\tau|)^{\ap/2s}<\e^{\ap}_{\vep}\bt/\|u_{\vep}\|_{L^{\iy}(\BR^n_I)}$, then by (4.15) we have the estimate
\begin{equation*}\begin{split}&|u_{\vep}(x,t)-u_{\vep}(y,\tau)|\\
&\quad\le\f{|u_{\vep}(x,t)-u_{\vep}(y,t)|}{|x-y|^{\ap}}|x-y|^{\ap}+\f{|u_{\vep}(y,t)-u_{\vep}(y,\tau)|}{|t-\tau|^{\ap/2s}}|t-\tau|^{\ap/2s}<2\bt,
\end{split}\end{equation*}
and thus we have that
\begin{equation}\vr(\bt)<2\bt.
\end{equation}
Let $\rho$ be the modulus of continuity of $u$ in the compact subset $(\overline{\Om^1})_{I_*}\s\Om_I$ of $\BR^n_{I_*}\s\Om_I$
defined by
\begin{equation*}\rho(\bt)=\sup_{\substack{(x,t),(y,\tau)\in(\overline{\Om^1})_{I_*}\s\Om_I\\(|x-y|^{2s}+|t-\tau|)^{1/2s}<\bt}}|u(x,t)-u(y,\tau)|.
\end{equation*}
Set $\vartheta_{\vep}=\vep+\vr(\vep)+\rho(\vep)$. By (4.15) and the continuity of $u$ in $\BR^n_{I_*}\s\Om_I$, we see that
\begin{equation}\lim_{\vep\downarrow 0}\vartheta_{\vep}=0.
\end{equation}

Furthermore, we have that
\begin{equation}u_{\vep}+\vartheta_{\vep}>u\,\,\text{ in $\BR^n_J\s\Om_I$.}
\end{equation}
Indeed, let us take any $(x,t)\in\BR^n_{I_*}\s\Om_I$. If $(x,t)\in[(\BR^n\s\Om_{\vep})\times I]\cup(\BR^n\times\{-T\})$, then we see that $u_{\vep}(x,t)=u(x,t)$, and so (4.18) works well. If $(x,t)\in(\Om_{\vep}\s\Om)\times I$, then there is some $y\in\pa\Om_{\vep}\subset\BR^n\s\Om_{\vep}$ such that $|x-y|\le\vep$. Thus we have that $u_{\vep}(y,t)=u(y,t)$, $|u_{\vep}(x,t)-u_{\vep}(y,t)|\le\vr(\vep)$ and $|u(x,t)-u(y,t)|\le\rho(\vep)$. Hence we obtain that
\begin{equation*}\begin{split}
u_{\vep}(x,t)-u(x,t)&\ge u_{\vep}(y,t)-u(y,t)-(\vr(\vep)+\rho(\vep))\\
&=-(\vr(\vep)+\rho(\vep))=-\vartheta_{\vep}+\vep>-\vartheta_{\vep},
\end{split}\end{equation*}
which gives (4.18).

If we set $w_{\vep}=u_{\vep}+\vartheta_{\vep}-u$, then we show that
\begin{equation*}w_{\vep}\ge 0\,\,\text{ in $\BR^n_I$.}
\end{equation*}
Indeed, $w_{\vep}$ is a weak supersolution of the nonlocal equation $\rL_K w_{\vep}+\pa_t w_{\vep}=0$ in $\Om_I$ with boundary condition (4.18), and thus it follows from Corollary 3.2.

Since $u\le u_{\vep}+\vartheta_{\vep}$ in $\BR^n_I$, by (4.11) we have that
\begin{equation}\begin{split}\e+u(x_0,t_0)&\le u(x_k,t_k)\le u_{\vep}(x_k,t_k)+\vartheta_{\vep} \\
&\le u_{\vep}(x_0,t_0)+2\vartheta_{\vep}
\end{split}\end{equation}
for a fixed $\vep\in(0,1]$ and a sufficiently large $k$.
By (4.13), we see that there is some $(y_{\vep},\tau_{\vep})\in\BR^n_{I_*}\s\Om^{\vep}_I$ such that
$(|y_{\vep}-x_0|^{2s}+|\tau_{\vep}-t_0|)^{1/2s}\le\vep.$
Since $u_{\vep}(y_{\vep},t_{\vep})=u(y_{\vep},t_{\vep})$, it follows from (4.19) that
\begin{equation}\e+u(x_0,t_0)\le u(y_{\vep},\tau_{\vep})+u_{\vep}(x_0,t_0)-u_{\vep}(y_{\vep},\tau_{\vep})+2\vartheta_{\vep}\le  u(y_{\vep},\tau_{\vep})+3\vartheta_{\vep}
\end{equation}
for all sufficiently small $\vep\in(0,1]$. Since $(y_{\vep},t_{\vep}),(x_0,t_0)\in(\overline{\Om^1})_{I_*}\s\Om_I$ and $\Om^{\vep}_I\subset\Om^1_I$ for $\vep\in(0,1]$, by the continuity of $u$ in $\BR^n_{I_*}\s\Om_I$ we have that
\begin{equation}u(y_{\vep},t_{\vep})-u(x_0,t_0)\le\rho(\vep)\le\vartheta_{\vep}.
\end{equation}
Thus by (4.20) and (4.21) we obtain that
$$\e+u(x_0,t_0)\le u(x_0,t_0)+4\vartheta_{\vep},$$
and so $\e\le 4\vartheta_{\vep}$. Taking $\vep\downarrow 0$, we have that $\e\le 0$ by (4.15), which gives a contradiction. Therefore we conclude that $u\in C(\BR^n_I)$. \qed

\begin{thm} If $u\in H^1(I;\rX_g(\Om))$ is a weak solution of the nonlocal parabolic equation ${\bf NP}_{\Om_I}(0,g,g)$ for $g\in C(\BR^n_{I_*})\cap L^{\iy}(\BR^n_I)$,
then $u$ is a viscosity solution of ${\,\bf NP}_{\Om_I}(0,g,g)$.
\end{thm}

\pf First, we show that any weak subsolution $u$ of the nonlocal equation ${\bf NP}_{\Om_I}(0,g,g)$ for $g\in C(\BR^n_{I_*})\cap L^{\iy}(\BR^n_I)$ is its viscosity subsolution. Take any $(x_0,t_0)\in\Om_I$. For a contrapositive proof, by continuity property we may assume that there are some $r>0$ with $Q_r(x_0,t_0)\subset\Om_I$ and $p\in \cP(\BR^{n+1})$ such that $u(x_0,t_0)=p(x_0,t_0)$, $u(x,t)<p(x,t)$ and
\begin{equation}\rL_K v(x,t)+\pa_t p(x,t)>0
\end{equation}
for any $(x,t)\in Q_r(x_0,t_0)\s\{(x_0,t_0)\}$,  where $v=p\mathbbm{1}_{Q_r(x_0,t_0)}+u\mathbbm{1}_{\BR^n\s Q_r(x_0,t_0)}$.
Then we see that $v\in C^2(Q_r(x_0,t_0))\cap L^{\iy}(\BR^n_I)$ by Theorem 4.1. Take any nonnegative testing function $\phi\in C^{\iy}_c(Q_r(x_0,t_0))$.
Then we have the estimate
\begin{equation}\begin{split}
&\int_{\BR^n}(\rL_K\fv(t))\phi\,dx \\
&\qquad=\int_{B_r(x_0)}\biggl(\int_{|y|<r-|x-x_0|}+\int_{|y|\ge r-|x-x_0|}\biggr)\mu_t(u,x,y)K(y)\,dy\,\phi(x,t)\,dx \\
&\qquad\le \f{c\|\n p\|_{L^{\iy}(Q_r(x_0,t_0))}}{1-s}\int_{B_r(x_0)}(r-|x-x_0|)^{2(1-s)}\phi(x,t)\,dx \\
&\qquad\qquad+\f{c\|v\|_{L^{\iy}(\BR^n_I)}}{s}\int_{B_r(x_0)}\f{\phi(x,t)}{(r-|x-x_0|)^{2s}}\,dx<\iy,
\end{split}\end{equation}
because $\overline{\{x\in B_r(x_0):\phi(x,t)\neq 0\}}\subset B_r(x_0)$ for any $t\in I_{r,s}(t_0)$. This implies that the integral in (4.23) is well-defined. Since $\fv(t)-\fu(t)=v(\cdot,t)-u(\cdot,t)\in\rX_0(B_r(x_0))$, $\fu(t)-\fg(t)=u(\cdot,t)-g(\cdot,t)\in\rX_0(\Om)$ and $\fg(t)=g(\cdot,t)\in H^s(\BR^n)$ for all $t\in I_{r,s}(t_0)$, it thus follows from Fubini's theorem, the change of variables, Lemma 5.1 and 5.2 \cite{CK} that
\begin{equation}\begin{split}
\int_{\BR^n}(\rL_K\fv(t))\phi\,dx&=\int_{\BR^n}\rL_K (\fv(t)-\fu(t))\,\phi\,dx \\ &\quad+\int_{\BR^n}\rL_K (\fu(t)-\fg(t))\,\phi\,dx+\int_{\BR^n}(\rL_K\fg(t))\phi\,dx \\
&=\la\fv(t)-\fu(t),\phi\ra_K+\la\fu(t)-\fg(t),\phi\ra_K+\la\fg(t),\phi\ra_K \\
&=\la\fv(t),\phi\ra_K.
\end{split}\end{equation}
Thus, by (4.22) and (4.24), we have that
\begin{equation*}\begin{split}\la\fv(t),\phi\ra_K+\int_{\BR^n}\fv'(t)\,\phi\,dx&=\la\fv(t),\phi\ra_K+\int_{\BR^n}(\pa_t v)\phi\,dx \\
&=\int_{B_r(x_0)}(\rL_K\fv(t)+\pa_t p)\phi\,dx\ge 0
\end{split}\end{equation*} for all $t\in I_{r,s}(t_0)$.
Thus $v$ is its weak supersolution on $Q_r(x_0,t_0)$, and so is $v+m$ on $Q_r(x_0,t_0)$ where
\begin{equation*}m=\inf_{\pa_p Q_r(x_0,t_0)}(u-p)<0.
\end{equation*}
By comparison principle (Corollary 3.3) on $Q_r(x_0,t_0)$, we have that $u\le v+m$ on $Q_r(x_0,t_0)$.
This gives a contradiction, because $u(x_0,t_0)\le v(x_0,t_0)+m<p(x_0,t_0)$.

Similarly, we can show that any weak subsolution of the nonlocal equation ${\bf NP}_{\Om_I}(0,g,g)$ for $g\in C(\BR^n_{I_*})\cap L^{\iy}(\BR^n_I)$ is its viscosity subsolution. Therefore we complete the proof. \qed

\begin{lemma} If $u\in H^1(I;\rX_g(\Om))$ is a weak solution of the nonlocal parabolic equation ${\bf NP}_{\Om_I}(0,g,g)$ for $g\in C(\BR^n_{I_*})\cap L^{\iy}(\BR^n_I)$,
then there is a universal constant $C>0$ such that
\begin{equation*}[u]_{C^{0,1}_x}(Q_r(x,t))\le\f{C}{(R-r)^{n+2s}}\|u\|_{L^{\iy}(\BR^n_I)}
\end{equation*} for any $r\in(0,R)$, where $Q_R(x,t)\subset\Om_I$.
\end{lemma}

\rk We apply Theorem 3.4 \cite{KL} and Theorem 5.2 \cite{KL1} in this proof. Looking over its proof scrupulously, we easily see that the H$\ddot{\rm o}$lder estimate holds for all $s\in(0,1)$ as follows; there are universal constants $C>0$ and $\ap\in(0,1)$ such that
\begin{equation}\|u\|_{C^{\ap}(Q_r(x,t))}\le C\bigl(\,\|u\|_{C(Q_R(x,t))}+\f{1}{(R-r)^{n+2s}}\|u\|_{L^{\iy}(\BR^n_I)}\bigr).
\end{equation}  for any $r\in(0,R)$ and $Q_R(x,t)\subset\Om_I$.

\,\,\pf Take any $r\in(0,R)$. For $k\in\BN\cup\{0\}$, set $r_k=r+2^{-k}(R-r)$. Then we see that $r_0=R>r_1>r_2>\cdots>r$ and $R-r_k\ge(R-r)/2$ for all $k\in\BN$.

Since $u$ is a viscosity solution by Theorem 4.3, by (4.25) there are universal constants $C>0$ and $\ap\in(0,1)$ such that
\begin{equation}\|u\|_{C^{\ap}(Q_{r_1}(x,t))}\le C\bigl(\,\|u\|_{C(Q_R(x,t))}+\f{1}{(R-r)^{n+2s}}\|u\|_{L^{\iy}(\BR^n_I)}\bigr).
\end{equation}
If we take $v=u\mathbbm{1}_{Q_R(x,t)}$ with $Q_R(x,t)\subset\Om_I$ as in the proof of Theorem 3.4 \cite{KL}, we obtain the estimate (4.26). Hence the required result follows from the standard telescopic sum argument \cite{CC}.
\qed

\section{Nonlocal weak Harnack inequality}

In this section, we shall prove nonlocal weak Harnack inequalities
with nonlocal parabolic tail term for weak subsolutions of the
nonlocal parabolic equation (1.3). This result plays a crucial role
in establishing nonlocal parabolic Harnack inequality for weak
subsolutions of the nonlocal parabolic equation (1.3).

To do this, first of all, we need the following {\it nonlocal Caccioppoli type inequality. }

\begin{thm} Let $\e\in C_c^{\iy}(t_0-r^{2s},\iy]$ be a function such that $0\le\e\le 1$ and $0\le\e'\le c/r^{2s}$ in $\BR$. If $u\in H^1(I;\rX^-_g(\Om))$ is a weak subsolution of the nonlocal parabolic equation {\bf NPE}$(0,g,g)$ given in $(1.3)$ and $Q^0_{2r}\subset\Om_I$ where $g\in H^s_T(\BR^n)$ and $f\le 0$ in $\Om_I$, and $w=(u-M)_{\pm}$ for $M\in\BR$,  then for any nonnegative $\zt\in C_c^{\iy}(B^0_r)$ we have the following estimate
\begin{equation*}\begin{split}
&\sup_{t\in I_{r,s}(t_0)}\e^2(t)\|w(\cdot,t)\zt\|^2_{L^2(B^0_r)}\\
&\quad+\int_{t_0-r^{2s}}^{t_0}\e^2(t)\iint_{B^0_r\times B^0_r}(\zt(x)w(x,t)-\zt(y)w(y,t))^2\,d_K(x,y)\,dt\\
&\le 2\int_{t_0-r^{2s}}^{t_0}\e(t)\e'(t)\|w(\cdot,t)\zt\|^2_{L^2(B^0_r)}\,dt\\
&\quad+\int_{t_0-r^{2s}}^{t_0}\e^2(t)\iint_{B^0_r\times B^0_r}[w(x,t)\vee w(y,t)]^2(\zt(x)-\zt(y))^2\,d_K(x,y)\,dt\\
&\quad+2\biggl(\,\sup_{(x,t)\in\supp(\zt)\times I_{r,s}(t_0)}\int_{\BR^n\s B^0_r}w(y,t)\,K(x-y)\,dy\biggr)\|w\zt^2\|_{L^1(Q^0_r)}.
\end{split}\end{equation*}
\end{thm}

\pf For simplicity of the proof, we may assume that
$(x_0,t_0)=(0,0)$. Let $w=(u-M)_+$ and take any $\zt\in
C_c^{\iy}(B_r)$. We use $\vp=w\zt^2\e^2$ as a testing function in
the weak formulation of the equation. Then we have that
\begin{equation*}\begin{split}
\int_{-r^{2s}}^\tau\int_{B_r}(\pa_t u)\vp\,dx\,dt+\rI(u,\vp)\le 0
\end{split}\end{equation*}
for any $\tau\in I_{r,s}(t_0)$, where the bilinear operator is given by
$$\rI(u,\vp)=\int_{-r^{2s}}^\tau\e^2(t)
\iint_{\Om\times\Om}(u(x,t)-u(y,t))(\vp(x,t)-\vp(y,t))\,d_K(x,y)\,dt.$$
The first term in the left-hand side of the above inequality can be evaluated by
\begin{equation}\begin{split}
&\int_{-r^{2s}}^\tau\int_{B_r}(\pa_t u)\vp\,dx\,dt=\int_{-r^{2s}}^\tau
\int_{B_r}\bigl[\pa_t(w^2\zt^2\e^2)-2 w^2\zt^2\e\e_t\bigr]dx\,dt\\
&=\e^2(\tau)\int_{B_r}w^2(x,\tau)\zt^2(x)\,dx-2\int_{-r^{2s}}^\tau\e(t)\e'(t)\int_{B_r}w^2(x,t)\zt^2(x)\,dx\,dt.
\end{split}\end{equation}
We next split $\rI(u,\vp)$ into two parts as follows;
\begin{equation}\begin{split}\rI(u,\vp)&=\int_{-r^{2s}}^\tau\e^2(t)\iint_{B_r\times B_r}(u(x,t)-u(y,t))(\vp(x,t)-\vp(y,t))d_K(x,y)dt\\
&+2\int_{-r^{2s}}^\tau\e^2(t)\int_{\BR^n\s B_r}\int_{B_r}(u(x,t)-u(y,t))\,\vp(x,t)\,d_K(x,y)\,dt\\
&:=I_1+2 I_2.
\end{split}\end{equation}
For the estimate of $I_1$, we first observe that
\begin{equation}(u(x,t)-u(y,t))(\vp(x,t)-\vp(y,t))\ge(w(x,t)-w(y,t))(\vp(x,t)-\vp(y,t))
\end{equation}
whenever $(x,t),(y,t)\in Q_r$; indeed, it can easily be checked by considering three possible cases (i) $u(x,t), u(y,t)>M$, (ii) $u(x,t)>M$, $u(y,t)\le M$, and (iii) $u(x,t)\le M$, $u(y,t)>M$.
For the estimate of $I_2$, we note that
\begin{equation*}\begin{split}
(u(x,t)-u(y,t))\vp(x,t)&\ge-(u(y,t)-u(x,t))_+(u(x,t)-M)_+\zt^2(x)\e^2(t)\\
&\ge-(u(y,t)-M)_+(u(x,t)-M)_+\zt^2(x)\e^2(t)\\
&=-w(y,t)w(x,t)\zt^2(x)\e^2(t)
\end{split}\end{equation*}
and thus we have that
\begin{equation}\begin{split}
I_2&\ge-\int_{-r^{2s}}^\tau\e^2(t)\int_{\BR^n\s B_r}\int_{B_r}w(y,t)w(x,t)\zt^2(x)\,d_K(x,y)\,dt\\
&\ge-\biggl(\,\sup_{(x,t)\in\supp(\zt)\times I_{r,s}}w(y,t)\,K(x-y)\,dy\biggr)\iint_{Q_r}w(x,t)\zt^2(x)\,dx\,dt.
\end{split}\end{equation}
Since the following equality
\begin{equation*}\begin{split}(w(x,t)-w(y,t))(\vp(x,t)-\vp(y,t))&=\e^2(t)(\zt(x)w(x,t)-\zt(y)w(y,t))^2\\
&\,\,\,-\e^2(t)w(x,t)w(y,t)(\zt(x)-\zt(y))^2
\end{split}\end{equation*}
is always true, it follows from (5.2), (5.3) and (5.4) that
\begin{equation}\begin{split}
&\rI(u,\vp)\ge\int_{-r^{2s}}^\tau\e^2(t)\iint_{B_r\times B_r}(\zt(x)w(x,t)-\zt(y)w(y,t))^2\,d_K(x,y)\,dt\\
&\,\,\,\,-\int_{-r^{2s}}^\tau\e^2(t)\iint_{B_r\times B_r}[w(x,t)\vee w(y,t)]^2(\zt(x)-\zt(y))^2\,d_K(x,y)\,dt\\
&\,\,\,\,-2\biggl(\,\sup_{(x,t)\in\supp(\zt)\times I_{r,s}}\int_{\BR^n\s B_r}w(y,t)\,K(x-y)\,dy\biggr)\iint_{Q_r}w(x,t)\zt^2(x)\,dx\,dt.
\end{split}\end{equation}
Hence the required inequality can be obtained from (5.1) and (5.5). \qed

\,\,Next, we obtain nonlocal weak Harnack inequality with nonlocal parabolic tail term for weak subsolutions of the nonlocal parabolic equation (1.3) in the following theorems.

\begin{thm} If $u\in H^1(I;\rX^-_g(\Om))$ is a weak subsolution of the nonlocal parabolic equation ${\bf NP}_{\Om_I}(f,g,g)$ given in $(1.3)$ and $Q^0_{2r}\subset\Om_I$ where $g\in H^s_T(\BR^n)$ and $f\le 0$ in $\Om_I$, then there is a constant $C_0=C_0(n,s,\ld,\Ld)>0$ such that
\begin{equation*}
\sup_{Q^0_r}u\le\dt\,\cT_r(u^+;(x_0,t_0))+C_0\,\dt^{-\f{\ap n}{4s}}\biggl(\displaystyle\f{1}{|\cQ^0_{2r}|}\iint_{Q^0_{2r}}[u^+]^2\,dx\,dt\biggr)^{\f{1}{2}}
\end{equation*} for any $\dt\in(0,1]$, where $\ap=1+\f{2s}{n}$.
\end{thm}

\pf Let $w=(u-M)_{\pm}$ for $M\in\BR$.
Using symmetry and the elementary inequlity $$[w(x,t)\vee w(y,t)]^2\le 2 w^2(x,t)+2 w^2(y,t),$$ by Theorem 5.1 and the mean value theorem we have that
\begin{equation}\begin{split}
&\sup_{t\in I_{r,s}(t_0)}\e^2(t)\|w(\cdot,t)\zt\|^2_{L^2(B^0_r)}\\
&\qquad+\int_{t_0-r^{2s}}^{t_0}\e^2(t)\iint_{B^0_r\times B^0_r}(\zt(x)w(x,t)-\zt(y)w(y,t))^2\,d_K(x,y)\,dt\\
&\le c\bigl(\,\|\e_t\|_{L^{\iy}(\BR)}+\|\n\zt\|^2_{L^{\iy}(B^0_r)}\bigr)\|w\|^2_{L^2(Q^0_r)}+\cA(w,\zt,t_0,r,s)\,\|w\|_{L^1(Q^0_r)}                      \\
\end{split}\end{equation} where $\cA(w,\zt,t_0,r,s)$ is the value given by
$$\cA(w,\zt,t_0,r,s)=2\sup_{(x,t)\in \supp(\zt)\times I_{r,s}(t_0)}\int_{\BR^n\s B^0_r}w(y,t)\,K(x-y)\,dy.$$
Applying a well-known parabolic version of the fractional Sobolev inequality to (5.6) and observing $|\cQ^0_r|/|Q^0_r|=2-\sm$, we obtain that
\begin{equation}\begin{split}
&\biggl(\f{1}{|\cQ^0_r|}\iint_{Q^0_r}|w\zt|^{2\ap}\,dx\,dt\biggr)^{\f{1}{\ap}}\\
&\le C_0 r^{2s}\bigl(\,\|\e'\|_{L^{\iy}(\BR)}+\|\n\zt\|^2_{L^{\iy}(B^0_r)}r^{2-2s}+r^{-2s}\bigr)\biggl(\f{1}{|\cQ^0_r|}\iint_{Q^0_r}|w|^2\,dx\,dt\biggr)\\
&\qquad\qquad+C_0 r^{2s}\cA(w,\zt,t_0,r,s)\,\f{1}{|\cQ^0_r|}\iint_{Q^0_r}w\,dx\,dt
\end{split}\end{equation} where $\ap=1+\f{2s}{n}$.
For $k=0,1,2,\cdots$, we set
$$r_k=(1+2^{-k})r,\,\,r^*_k=\f{r_k+r_{k+1}}{2},\,\,M_k=M+(1-2^{-k})M_*,\,\,M^*_k=\f{M_k+M_{k+1}}{2},$$
$w_k=(u-M_k)_+$ and $w^*_k=(u-M^*_k)_+$ for a constant $M_*>0$ to be determined later.
In (5.6), for $k=0,1,\cdots,$ we choose a function $\e_k\in C_c^{\iy}(t_0-(r_k^*)^{2s},\iy)$ with $\e_k|_{[t_0-r^{2s}_{k+1},t_0)}\equiv 1$ such that $0\le\e_k\le 1$ and $0\le\e_k'\le 2^{(k+2)2s}r^{-2s}$ in $\BR$, and a function $\zt_k\in C_c^{\iy}(B^0_{r^*_k})$ with $\zt_k|_{B^0_{r_{k+1}}}\equiv 1$ such that $0\le\zt_k\le 1$ and $|\n\zt_k|\le 2^{k+2}/r$ in $\BR^n$.
For $k=0,1,2,\cdots$, we set
$$N_k=\biggl(\f{1}{|\cQ^0_{r_k}|}\iint_{Q^0_{r_k}}|w_k|^2\,dx\,dt\biggr)^{\f{1}{2}}.$$
Since $w_k^*\ge w_{k+1}$ and $w_k^*(x,t)\ge M_{k+1}-M_k^*=2^{-k-2}M_*$ whenever $u(x,t)\ge M_{k+1}$, we then have that
\begin{equation}\begin{split}
N_{k+1}&\le c\biggl(\f{1}{|\cQ^0_{r_k}|}\iint_{Q^0_{r_{k+1}}}\f{w^2_{k+1}(w_k^*)^{2(\ap-1)}}{(M_{k+1}-M_k^*)^{2(\ap-1)}}\,dx\,dt\biggr)^{\f{1}{2}} \\
&\le c\biggl(\f{2^k}{M_*}\biggr)^{\ap-1}\biggl(\f{1}{|\cQ^0_{r_k}|}\iint_{Q^0_{r_k}}|w^*_k\,\zt_k|^{2\ap}\,dx\,dt\biggr)^{\f{1}{2}}
\end{split}\end{equation}
Since $\zt_k\in C_c^{\iy}(B^0_{r_k})$, $w_k^*\le w_0$ for all $k$ and
$$|y-x|\ge|y-x_0|-|x-x_0|\ge\biggl(1-\f{r_k^*}{r_k}\biggr)|y-x_0|\ge 2^{-k-2}|y-x_0|$$
for any $x\in B^0_{r^*_k}$ and $y\in\BR^n\s B^0_{r_k}$, we easily obtain that
\begin{equation}
\cA(w_k^*,\zt_k,t_0,r_k,s)\le c\,2^{k(n+2s)}\,r^{-2s}\,\cT_r(w_0;(x_0,t_0)).
\end{equation}
Since $0\le w_k^*\le w_k$ and $w_k(x,t)\ge M_k^*-M_k=2^{-k-2}M_*$ if $u(x,t)\ge M_k^*$, it follows from (5.7), (5.8) and (5.9) that
\begin{equation*}\begin{split}
&\biggl(\f{2^k}{M_*}\biggr)^{-\f{2(\ap-1)}{\ap}}N^{2/\ap}_{k+1}
\le\f{c\,2^{2k}}{|\cQ^0_{r_k}|}\iint_{Q^0_{r_k}}|w^*_k|^2\,dx\,dt \\
&\qquad\qquad\qquad+c\,2^{k(n+2s)}\biggl(\f{r_k}{r}\biggr)^{2s}\,\cT_r(w_0;(x_0,t_0))\,\f{1}{|\cQ^0_{r_k}|}\iint_{Q^0_{r_k}}w^*_k\,dx\,dt\\
&\qquad\quad\le c\,2^{2k}N_k^2+c\,2^{k(n+2s)}\,\cT_r(w_0;(x_0,t_0))\f{1}{|\cQ^0_{r_k}|}\iint_{Q^0_{r_k}}\f{w^*_k w_k}{M^*_k-M_k}\,dx\,dt\\
&\quad\qquad\le c\biggl(2^{2k}+ \f{\,2^{k(n+2s+1)}}{M_*}\,\cT_r(w_0;(x_0,t_0))\biggr)\,N_k^2.
\end{split}\end{equation*}
Taking $M^*$ in the above so that
$$M_*\ge \dt\,\cT_r(w_0;(x_0,t_0))\,\,\text{ for $\dt\in(0,1]$,}$$
we obtain that
\begin{equation}\f{N_{k+1}}{M_*}\le d_0\,a^k\,\biggl(\f{N_k}{M_*}\biggr)^{1+\f{2s}{n}}
\end{equation}
where $d_0=c^{\f{\ap}{2}}\dt^{-\f{\ap}{2}}$ and $a=2^{\f{\ap}{2}(n+2s+1)+\f{2s}{n}}$. If $N_0\le d_0^{-\f{n}{2s}}a^{-\f{n^2}{4 s^2}}M_*,$ then we set
$$M_*=\dt\,\cT_r(w_0;(x_0,t_0))+c_0\,\dt^{-\f{\ap n}{4s}}a^{\f{n^2}{4 s^2}}N_0$$ where $c_0=c^{\f{\ap n}{4s}}$.
By Lemma 2.2, we conclude that
\begin{equation*}\begin{split}
\sup_{Q^0_r}u&\le M+M_*\\
&\le M+\dt\,\cT_r(w_0;(x_0,t_0))+c_0^{\f{\ap n}{4s}}\dt^{-\f{\ap n}{4s}}a^{\f{n^2}{4 s^2}}\biggl(\displaystyle\f{1}{|\cQ^0_{2r}|}\iint_{Q^0_{2r}}(u-M)_+^2\,dx\,dt\biggr)^{\f{1}{2}}.
\end{split}\end{equation*}
Hence, choosing $M=0$ in the above estimate, we obtain the required result. \qed

\,\,\,The third one is a lemma which furnishes a precise relation between the nonlocal parabolic tails of the positive and negative part of the weak subsolutions.

\begin{lemma} If $u\in H^1(I;\rX^-_g(\Om))$ is a weak subsolution of the nonlocal parabolic equation ${\bf NP}_{\Om_I}(f,g,g)$ in $(1.3)$ such that $u\ge 0$ in $Q^0_R\subset\Om_I$ where $g\in C(\BR^n_{I_*})\cap L^{\iy}(\BR^n_I)$ and $f\le 0$ in $\Om_I$, then we have the estimate
\begin{equation*}\cT_r(u^+;(x_0,t_0))\lesssim\,\sup_{Q^0_r}u+\biggl(\f{r}{R}\biggr)^{2s}\cT_R(u^-;(x_0,t_0))
\end{equation*}
for any $r$ with $0<r<R$.
\end{lemma}

\pf  Without loss of generality, we may assume that $(x_0,t_0)=(0,0)$ and $\cT_r(u^+;(0,0))\neq 0$. Fix any $\vep>0$ with $\vep\le\cT_r(u^+;(0,0))/2$. Then it follows from the definition of supremum, the uniform continuity of $u^+$ on a big enough closed ball (via Theorem 4.2) and the Lebesgue's dominated convergence theorem (via Theorem 4.1) that
there are some $\tau\in(-r^{2s},0]$ and $\ep_0>0$ with $[\tau-\ep_0 r^{2s},\tau]\subset (-r^{2s},0]$
such that
\begin{equation}
r^{2s}\int_{\tau-\ep_0 r^{2s}}^{\tau}\int_{\BR^n\s B_r}\f{u^+(y,t)}{|y|^{n+2s}}\,dx
\ge \cT_r(u^+;(0,0))-\vep\ge\f{1}{2}\,\cT_r(u^+;(0,0)).
\end{equation}
Let $M=\sup_{Q_r}u$ and set $\vp(x,t)=(u(x,t)-2M)\zt^2(x)\e^2(t)$ where $\zt\in C_c^{\iy}(B_{3r/4})$ is a function satisfying that $\zt|_{B_{r/2}}\equiv 1$, $0\le\zt\le 1$ and $|\n\zt|\le c/r$ in $\BR^n$, and $\e\in C_c^{\iy}(-r^{2s},0]$ is a function such that $\e=1$ in $[\tau-\ep_0 r^{2s},\tau]$,  $0\le\e\le 1$ and $|\e'|\le c/r^{2s}$ in $\BR$. Then we have that
\begin{equation}\begin{split}
0&\ge-\iint_{\Om\times I}u\,\pa_t\vp\,dx\,dt\\
&+\int_{-r^{2s}}^0\iint_{B_r\times B_r}(u(x,t)-u(y,t))(\vp(x,t)-\vp(y,t))\,d_K(x,y)\,dt\\
&+2\int_{-r^{2s}}^0\e^2(t)\int_{B^c_r}\int_{B_r}(u(x,t)-u(y,t))(u(x,t)-2M)\zt^2(x)d_K(x,y)dt\\
&:=J_1+J_2+2 J_3.
\end{split}\end{equation}
Since the fact that $$-2M\le w(x,t):=u(x,t)-2M\le-M$$ for any $(x,t)\in Q_r$ and the following elementary equality
$$(\bt-\ap)(B^2\bt-A^2\ap)=(B\bt-A\ap)^2-\ap\bt(B-A)^2\text{ for any $\ap,\bt\in\BR$ and $A,B\ge 0$ }$$leads us to obtain the estimate
$$(w(x,t)-w(y,t))(w(x,t)\zt^2(x)-w(y,t)\zt^2(y))\ge-4\,M^2(\zt(x)-\zt(y))^2$$
for any $(x,t),(y,t)\in Q_r$,
it follows from simple calculation that
\begin{equation}\begin{split}
J_1+J_2&\ge \biggl[-\f{\e^2(t)}{2}\int_{B_r}u^2(x,t)\zt^2(x)\,dx\biggr]_{t=-r^{2s}}^{t=0}\\
&-2\int_{-r^{2s}}^0\e(t)\e'(t)\int_{B_r}\bigl[u(x,t)(u(x,t)-2M)-u^2(x,t)\bigr]dx\,dt\\
&-4 M^2\int_{-r^{2s}}^0\e^2(t)\iint_{B_r\times B_r}(\zt(x)-\zt(y))^2\,d_K(x,y)\,dt\\
&\gtrsim-M^2 r^{-2s}|Q_r|.
\end{split}\end{equation} The lower estimate on $J_3$ can be splitted as follows;
\begin{equation*}\begin{split}
J_3&\ge\int_{-r^{2s}}^0\e^2(t)\int_{\BR^n\s B_r}\int_{B_r}M(u(y,t)-M)_+\,\zt^2(x)\,d_K(x,y)\,dt\\
&-2M\int_{-r^{2s}}^0\e^2(t)\int_{E^t_M}\int_{B_r}(u(x,t)-u(y,t))_+\,\zt^2(x)\,d_K(x,y)\,dt\\
&:=J_3^1-J_3^2,
\end{split}\end{equation*} where $E^t_M=\{y\in\BR^n\s B_r:u(y,t)<M\}$.
Since $(u(y,t)-M)_+\ge u^+(y,t)-M$, by (5.11) the lower estimate on $J_3^1$ can here be obtained as
\begin{equation}\begin{split}
J_3^1&\ge d_2 M r^{-2s}|Q_r|\,\cT_r(u^+;(0,0))-d_3 M^2 r^{-2s}|Q_r|
\end{split}\end{equation}
with universal constants $d_2,d_3>0$. If $(x,t)\in Q_r$ and $y\in E^t_M$, then we observe that
\begin{equation*}\begin{split}(u(x,t)-u(y,t))_+&\le |u(x,t)-M|+|M-u(y,t)|\\
&\le M+(M+u^-(y,t)-u^+(y,t))\\
&\le 2M+u^-(y,t)
\end{split}\end{equation*}
because $u^+(y,t)< M+u^-(y,t)$ for any $y\in E^t_M$.
Since $u^-(y,t)=0$ for all $(y,t)\in Q_R$, the upper estimate on $J_3^2$ can thus be achieved by
\begin{equation}\begin{split}
J_3^2&\le 2 M^2\int_{-r^{2s}}^0\e^2(t)\int_{\BR^n\s B_r}\int_{B_r}\zt^2(x)\,d_K(x,y)\\
&+2M\int_{-r^{2s}}^0\e^2(t)\int_{B_R\s B_r}\int_{B_r}(2M+u^-(y,t))\zt^2(x)\,d_K(x,y)\\
&+2M\int_{-r^{2s}}^0\e^2(t)\int_{\BR^n\s B_R}\int_{B_r}(2M+u^-(y,t))\zt^2(x)\,d_K(x,y)\\
&\le d_4 M^2 r^{-2s}|Q_r|+d_5 M R^{-2s}|Q_r|\,\cT_R(u^-;(0,0))
\end{split}\end{equation}
with universal constants $d_4,d_5>0$. Thus, by (5.14) and (5.15), we have that
\begin{equation}\begin{split}
J_3&\ge -d M^2 r^{-2s}|Q_r|-d M R^{-2s}|Q_r|\,\cT_R(u^-;(0,0))\\
&\quad+e M r^{-2s} |Q_r|\,\cT_r(u^+;(0,0))
\end{split}\end{equation}
where $d,e>0$ are some universal constants depending only on $n,s,\ld$ and $\Ld$.
Hence the estimates (5.12), (5.13) and (5.16) give the required estimate. \qed

\,\,\,Next we obtain a weak Harnack inequality for nonnegative weak subsolutions of the nonlocal parabolic equation (1.3) by employing Theorem 5.2 and Lemma 5.3. It is interesting that this estimate no longer depends on the nonlocal parabolic tail term, but its proof is pretty simple.

\begin{thm} If $u\in H^1(I;\rX^-_g(\Om))$ is a nonnegative weak subsolution of the nonlocal parabolic equation ${\bf NP}_{\Om_I}(f,g,g)$ given in $(1.3)$ where $g\in C(\BR^n_{I_*})\cap L^{\iy}(\BR^n_I)$ and $f\le 0$ in $\Om_I$, then we have the estimate
$$\sup_{Q^0_r} u\le C\biggl(\f{1}{|Q^0_{2r}|}\iint_{Q^0_{2r}}u^2(x,t)\,dx\,dt\biggr)^{1/2}$$
for any $r>0$ with $Q^0_{2r}\subset\Om_I$.
\end{thm}

\pf We choose some $\dt\in(0,1]$ so that $1-\dt d_0>0$ and take any $r>0$ with $Q^0_{2r}\subset\Om_I$. Then it follows from Theorem 5.2 and Lemma 5.3 that
\begin{equation*}\begin{split}
\sup_{Q^0_r}u&\le\dt\,d_0\biggl[\,\sup_{Q^0_r}u+\cT_r(u^-;(x_0,t_0))\biggr]\\
&\qquad+C_0\,\dt^{-\f{\ap n}{4s}}\biggl(\f{1}{|Q^0_{2r}|}\iint_{Q^0_{2r}}u^2(x,t)\,dx\,dt\biggr)^{\f{1}{2}}
\end{split}\end{equation*}
Since $\cT_r(u^-;(x_0,t_0))=0$, we can easily derive the required result by taking
$$C=\f{C_0\,\dt^{-\f{\ap n}{4s}}}{1-\dt d_0}.$$
Hence we complete the proof. \qed

\section{Parabolic fractional Poincar\'e inequality}

Let $n\ge 1$, $p\ge 1$, $s\in(0,1)$ and $sp<n$. For a ball $B\subset\BR^n$, let $u_B$ denote the average of $u\in W^{s,p}(B)$ over $B$, i.e.
$$u_B=\f{1}{|B|}\int_B u(y)\,dy.$$ Then it is well-known [BBM, MS] that
\begin{equation}\|u-u_B\|^p_{L^p(B)}\le\f{c_{n,p}(1-s)|B|^{\f{sp}{n}}}{(n-sp)^{p-1}}\,\|u\|^p_{W^{s,p}(B)}
\end{equation}
with a universal constant $c_{n,p}>0$ depending only on $n$ and $p$,
which is called the {\it fractional Poincar\'e inequality}. We note
that this estimate no longer depends on partial differential
equations. However, the inequality as (6.1) in the parabolic sense
is not available for a general function $u(x,t)\in
L^2(I;\rX_0(\Om))$.

If $u\in L^{\iy}(\BR^n_I)$, then we have that $\cT_R(u;(x_0,t_0))\le\|u\|_{L^{\iy}(\BR^n_I)}$. So, if we could assume that $ \|u\|_{L^{\iy}(\BR^n_I)}=1$, we see that
$$\kappa:=\f{1}{2}\biggl(\f{r}{R}\biggr)^{2s}\cT_R(u^-;(x_0,t_0))<1\,\,\text{ for any $r\in(0,R)$,} $$
where $R>0$ satisfies $Q^0_R\subset\Om_I$.
For $Q^0_r:=Q_r(x_0,t_0)$, we denote by  $$u_{Q_r^0}=\f{1}{|Q_r^0|}\iint_{Q_r^0}u(x,t)\,dx\,dt.$$
For $x_0\in\BR^n$ and $\kappa\in(0,1)$, we consider a radial function $\theta\in C_c^{\iy}(\BR^n)$ with values in $[0,\kappa]$ such that $\theta|_{B_{r}(x_0)}\equiv\sqrt\kappa$, $\theta|_{\BR^n\s{B_{2r}(x_0)}}\equiv 0$ and $|\n\theta|\le c\sqrt\kappa/r$ in $\BR^n$. For simplicity, we write $d\theta(x,t)=\theta(x)\,dx\,dt$, $d\theta(x)=\theta(x)\,dx$,
$$\theta(Q^0_r)=\iint_{Q^0_r}d\theta(x,t)\,\,\,\,\text{ and }\,\,\,\,\theta(B^0_r)=\int_{B^0_r}d\theta(x).$$
Also, we denote by
$$u_{Q^0_r}^{\theta}=\f{1}{\theta(Q^0_r)}\iint_{Q^0_r}u(x,t)\,d\theta(x,t)\text{ and }u_{B^0_r}^{\theta}(t)=\f{1}{\theta(B^0_r)}\int_{B^0_r}u(x,t)\,d\theta(x).$$

Now we establish a parabolic version of the fractional Poincar\'e inequality for weak solutions for nonlocal parabolic equation (1.3) as follows.

\begin{thm} Let $u\in H^1(I;\rX_g(\Om))$ be an weak solution of the nonlocal parabolic equation ${\bf NP}_{\Om_I}(0,g,g)$ such that $u\ge\kappa>0$ in $Q^0_R\subset\Om_I$ where $g\in C(\BR^n_{I_*})\cap L^{\iy}(\BR^n_I)$ and
$$\kappa=\f{1}{2}\biggl(\f{r}{R}\biggr)^{2s}\cT_R(u^-;(x_0,t_0)).$$ If $\,v(x,t)=\ds\ln\biggl(\f{d}{u(x,t)}\biggr)$ for $d>0$ and there is a $\nu\in(0,1)$ such that
\begin{equation}|Q^0_r\cap\{u\ge d+\kappa\}|\ge\nu\,|Q^0_r|
\end{equation} for some $r$ with $0<5r<R$,
then for each $\sm\in(0,\nu\wedge 2/5)$ there exists a constant $C_0=C_0(n,s,\ld,\Ld,\nu,\sm)>0$ as follows;

\,$(a)$ there is some $\tau_0\in[t_0-r^{2s},t_0-\sm r^{2s}]$ such that
\begin{equation}\sup_{b\in[t_0-\sm r^{2s},t_0)}\int_{\tau_0}^b\iint_{B_{2r}^0\times B_{2r}^0}\f{|v(x,t)-v(y,t)|^2}{|x-y|^{n+2s}}\,dx\,dy\,dt\le C_0\,r^n,
\end{equation}

\,$(b)$ we have the estimate
\begin{equation}\begin{split}
&\iint_{(Q^0_r)^+}|v(x,t)-v_{(Q^0_r)^+}|^2\,dx\,dt\\
&\qquad\le C_0\,r^{2s}\int_{t_0-\sm r^{2s}}^{t_0}\iint_{B^0_{2r}\times B^0_{2r}}\f{|v(x,t)-v(y,t)|^2}{|x-y|^{n+2s}}\,dx\,dy\,dt\\
&\quad\qquad+C_0\,|(Q^0_r)^+|
\end{split}\end{equation} where $(Q^0_r)^+:=B^0_r\times[t_0-\sm r^s,t_0)=B_r(x_0)\times[t_0-\sm r^s,t_0)$ for $r>0$.
\end{thm}

\rk The reason why $\sm\le 2/5$ follows from the inequality (7.42) below.

\,\,We now state three lemmas which is useful for the proof of Theorem 6.1. The first one is an elliptic version of {\it weighted Poincar\'e inequality} which was obtained by [FK] as follows.

\begin{lemma} If $u\in W^{s,2}(B^0_{2r})$ for $s\in(0,1)$ and $r>0$, then there is a constant $c=c(n,s,\Ld)>0$ such that
\begin{equation*}
\int_{B^0_{2r}}|u(x)-u^{\theta}_{B^0_{2r}}|^2\,\theta(x)\,dx\le c\,r^{2s}\iint_{B_{2r}^0\times B_{2r}^0}|u(x)-u(y)|^2\bigl(\theta(x)\wedge\theta(y)\bigr)d_K(x,y)
\end{equation*} where $d_K(x,y)=K(x-y)\,dx\,dy$.
\end{lemma}

The second one is the following lemma whose detailed proof can be found in \cite{L}.

\begin{lemma} If $u\in L^2((Q^0_r)^+)$ for $r>0$, then we have that
$$\iint_{(Q^0_r)^+}|u(x,t)-u_{(Q_r^0)^+}|^2\,dx\,dt\le 4\iint_{(Q_r^0)^+}|u(x,t)-h|^2\,dx\,dt$$
for any $h\in\BR$.
\end{lemma}

\begin{lemma} For $d>0$, the functions $r_{d,\ep}(t)=\ds\ln^+\biggl(\f{d}{t}\biggr)\bigl((d-t)_++\ep\bigr)^{-2}$ and $q_{d,\ep}(t)=\ds\ln^+\biggl(\f{d}{t}\biggr)(d-t)_+\bigl((d-t)_++\ep\bigr)^{-3}$ are decreasing on $(0,\iy)$ for $\ep\ge (2 e^{-\f{1}{2}}-1)d$ and are Lipschitz continuous in $(\kappa,\iy)$ for $\kappa\in(0,1)$.
\end{lemma}

\pf Since $r_{d,\ep}(t)=0$ on $t\ge d$, we have only to consider $0<t<d$. Then we have that
$$r'_{d,\ep}(t)=-(d+\ep-t)^{-3}\biggl(\f{d+\ep-t}{t}-2\ln\biggl(\f{d}{t}\biggr)\biggr):=-(d+\ep-t)^{-3}\,h(t).
$$
Since $h'(t)=t^{-2}\bigl(2t-(d+\ep)\bigr)$ and
$$h\biggl(\f{d+\ep}{2}\biggr)=1-2\ln\biggl(\f{2d}{d+\ep}\biggr)\ge 0$$ for $\ep\ge (2 e^{-\f{1}{2}}-1)d$, we see that $h(t)\ge 0$, and so $r'_{d,\ep}(t)\le 0$.
Note that $q_{d,\ep}(t)=r_{d,\ep}(t)p_{d,\ep}(t)$ where
\begin{equation}p_{d,\ep}(t)=\f{(d-t)_+}{(d-t)_++\ep}.
\end{equation}
Since it is easy to check that $p_{d,\ep}(t)$ is decreasing on $(0,\iy)$ for any $\ep>0$, we can easily conclude the first required result. Also the second required result can be easily  obtained. Hence we are done. \qed

\,\,\,{\bf [Proof of Theorem 6.1]}
Without loss of generality, we may assume that $(x_0,t_0)=(0,0)$. Take any $a,b\in I_{r,s}$ with $a<b$ and choose $r>0$ so that $B_{5r}\subset B_R$ where $Q_R\subset\Om_I$. Then we have two possible cases; either $r>1$ or $r\le 1$.
\begin{equation}\begin{cases}B_{4r+h_0}\subset B_R &\text{ for $ r>1$ } \\
r\le r^{n+2s} &\text{ for $r\le 1,$}
\end{cases}\end{equation}
where $h_0=r^{\f{1}{n+2s}}$.
We use the function $$\vp(x,t)=\f{\theta^2_0(x,t)}{\kappa\,u(x,t)}$$ as a testing function where $\theta_0(x,t)=\theta(x/2)h(x,t)$ and $$h(x,t)=\f{\sqrt\kappa\,(2-u(x,t))_+}{(2-u(x,t))_++\ep_0}$$ for $\ep_0=2(2 e^{-\f{1}{2}}-1)$. From (1.3), we obtain that
\begin{equation}\int_{\Om}\pa_t u(x,t)\,\vp(x,t)\,dx+\cI_t(u,\vp)=0
\end{equation} for $\aee$ $t\in I$,
where the bilinear operator $\cI_t$ is given by
$$\cI_t(u,\vp)=\iint_{\BR^n\times\BR^n}(u(x,t)-u(y,t))(\vp(x,t)-\vp(y,t))\,d_K(x,y).$$
Since $v(x,t)-v(y,t)=V(x,t)-V(y,t)$ and $v(x,t)-v_{Q^+_r}=V(x,t)-V_{Q^+_r}$ where $$V(x,t)=\ln\biggl(\f{d\,\|u\|_{L^{\iy}(\BR^n_I)}}{u(x,t)}\biggr),$$
without loss of generality, we may assume that $\|u\|_{L^{\iy}(\BR^n_I)}=1$; for, $u/\|u\|_{L^{\iy}(\BR^n_I)}$ is also an weak solution of the equation ${\bf NP}_{\Om_I}(0,g,g)$, and we use it in place of $u$.
Now the integral of the first term on $[a,b]\subset I_{r,s}$ in (6.7) can be estimated as
\begin{equation}\begin{split}
\int_a^b\int_{\Om}\pa_t u(x,t)\,\vp(x,t)\,dx\,dt&=-\int_a^b\int_{\Om}\pa_t\biggl[\f{1}{\kappa}\ln\biggl(\f{2}{u(x,t)}\biggr)\biggr]\,\theta^2_0(x,t)\,dx\,dt\\
&=-\int_{\Om}\int_a^b\pa_t v_0(x,t)h^2(x,t)\,dt\,\theta^2_1(x)\,dx\\
&=-\int_\Om\int_a^b\pa_t[v_0(x,t)h^2(x,t)]\,dt\,\theta^2_1(x)\,dx     \\
&\qquad+\int_\Om\int_a^b v_0(x,t)\,\pa_t h^2(x,t)\,dt\,\theta^2_1(x)\,dx     \\
&:=(w^{\theta^2_1}_{B_{4r}}(a)-w^{\theta^2_1}_{B_{4r}}(b))\theta^2_1(B_{4r})+\cI\\
\end{split}\end{equation}
where  $\theta_1(x)=\theta(x/2)$, $w(x,t)=v_0(x,t)h^2(x,t)$ for $v_0(x,t)=\f{1}{\kappa}\ln^+(2/u(x,t))$ and
$$\cI=\int_{\Om}\int_a^b v_0(x,t)\,\pa_t h^2(x,t)\,dt\,\theta^2_1(x)\,dx.$$
Here we observe that
\begin{equation}-\cI=\int_a^b\int_{\Om} \pa_t u(x,t)\,\psi(x,t)\,\theta^2_1(x)\,dx\,dt=-\int_a^b\cI_t(u,\psi\,\theta^2_1)\,dt
\end{equation}
where $$\psi(x,t)=\f{\ep_0\,(2-u(x,t))_+\,v_0(x,t)}{[(2-u(x,t))_++\ep_0]^3}.$$
Since $\ln^+\e\le(\e-1)_+$ for all $\e>0$, we have that
\begin{equation}\psi(x,t)\le\f{\ep_0(2-u(x,t))^2_+}{\,u(x,t)[(2-u(x,t))_++\ep_0]^3}\le\f{1}{\kappa}.
\end{equation}
We note that
\begin{equation*}\begin{split}&\bigl(u(x,t)-u(y,t)\bigr)\bigl(\psi(x,t)\theta^2_1(x)-\psi(y,t)\theta^2_1(y)\bigr)\\
&\qquad\qquad\qquad=\bigl(u(x,t)-u(y,t)\bigr)\bigl(\psi(x,t)-\psi(y,t)\bigr)\theta^2_1(x)\\
&\qquad\qquad\qquad\quad+\bigl(u(x,t)-u(y,t)\bigr)\bigl(\theta^2_1(x)-\theta^2_1(y)\bigr)\psi(y,t).
\end{split}\end{equation*}
By Lemma 6.4, we have that
\begin{equation}\cI_t(u,\psi)\le 0\,\,\text{ for any $t\in I$,}\end{equation} where
$$\cI_t(u,\psi)=\iint_{\BR^n\times\BR^n}\bigl(u(x,t)-u(y,t)\bigr)\bigl(\psi(x,t)-\psi(y,t)\bigr)\,\theta^2_1(x)\,d_K(x,y).$$
Indeed, assuming without loss of generality that $u(x,t)\ge u(y,t)$, it follows from Lemma 6.4 that $$(u(x,t)-u(y,t))(\psi(x,t)-\psi(y,t))\le 0$$ for any $x,y\in\Om$ and $t\in I$. Thus we conclude that $\cI_t(u,\psi)\le 0$.

Also, we can derive the estimate $\,\cJ_t(u,\theta^2_1)\le c\,r^{n-2s}$ where
\begin{equation*}\cJ_t(u,\theta_1)=\iint_{\BR^n\times\BR^n}\bigl(u(x,t)-u(y,t)\bigr)\bigl(\theta^2_1(x)-\theta^2_1(y)\bigr)\psi(y,t)\,d_K(x,y).
\end{equation*}
Indeed, by the mean value theorem, we have that
\begin{equation*}\begin{split}u(x,t)-u(y,t)&=\int_0^1\la \n u(x+\tau(y-x),t),y-x\ra\,d\tau, \\
\theta^2_1(x,t)-\theta^2_1(y,t)&=\int_0^1\theta\biggl(\f{x+\tau(y-x)}{2}\biggr)\biggl\la \n \theta\biggl(\f{x+\tau(y-x)}{2}\biggr),y-x\biggr\ra\,d\tau,
\end{split}\end{equation*}
and so it follows from  the definition of $\theta_1$, (6.6), (6.10) and Lemma 4.4 that
\begin{equation}\begin{split}\cJ_t(u,\theta^2_1)&\le\f{c}{r^2}\,\|u\|_{L^{\iy}(\BR^n_I)}\iint_{B_{5r}\times B_{5r}}\f{|x-y|^2}{|x-y|^{n+2s}}\,dx\,dy \\
&\,\,+2\|u\|_{L^{\iy}(\BR^n_I)}\iint_{(\BR^n\s B_{5r})\times B_{4r}}\f{1}{|x-y|^{n+2s}}\,dx\,dy\le c\,r^{n-2s},
\end{split}\end{equation} because $\|u\|_{L^{\iy}(\BR^n_I)}=1$.
Thus, by (6.9), (6.11) and (6.12), we obtain that
\begin{equation}\begin{split}-\cI&=-\int_a^b\cI_t(u,\psi\,\theta^2_1)\,dt=-\int_a^b\cI_t(u,\psi)\,dt-\int_a^b\cJ_t(u,\theta^2_1)\,dt \\
&\ge-\int_a^b\cJ_t(u,\theta^2_1)\,dt\ge-c\,r^n,
\end{split}\end{equation}
and thus it follows from (6.8) and (6.13) that
\begin{equation}(w^{\theta^2_1}_{B_{4r}}(a)-w^{\theta^2_1}_{B_{4r}}(b))\theta^2_1(B_{4r})+ c r^n\ge -\int_a^b\cI_t(u,\vp)\,dt.
\end{equation}
We note that
\begin{equation*}\begin{split}
\cI_t(u,\vp)&=2\iint_{(\BR^n\s B_{4r})\times B_{4r}}(u(x,t)-u(y,t))\,\vp(x,t)\,d_K(x,y)  \\
&\quad+\iint_{B_{4r}\times B_{4r}}(u(x,t)-u(y,t))(\vp(x,t)-\vp(y,t))\,d_K(x,y) \\
&:=I_t(u,\vp)+J_t(u,\vp).
\end{split}\end{equation*}
Since $h^2(y,t)\ge e/16$ for any $(y,t)\in B_{2r}\times[a,b]$, as in the proof of Lemma 1.3 \cite{DKP1} we have the following estimates

\begin{equation}\begin{split}
I_t(u,\vp)&\le c\,r^{n-2s}+c\,\kappa^{-1}r^n \int_{\BR^n\s B_R}\f{u^-(y,t)}{|y|^{n+2s}}\,dy,\\
J_t(u,\vp)&\le-\f{1}{c\,\kappa}\iint_{B_{4r}\times B_{4r}}\biggl[\ln\biggl(\f{d}{u(x,t)}\biggr)-\ln\biggl( \f{d}{u(y,t)}\biggr)\biggr]^2\theta^2_0(y,t)\,d_K(x,y)\\
&\qquad+c\iint_{B_{4r}\times B_{4r}}(\theta_0(x,t)-\theta_0(y,t))^2\,d_K(x,y)\\
&\le-\f{e}{16c}\iint_{B_{2r}\times B_{2r}}\f{|v(x,t)-v(y,t)|^2}{|x-y|^{n+2s}}\,dx\,dy+c\,J_t(\theta_0,\theta_0).
\end{split}\end{equation}
For the estimate of $I_t(\theta_0,\theta_0)$, we note that
\begin{equation}(\theta_0(x,t)-\theta_0(y,t))^2\le 2(\theta(x)-\theta(y))^2 h^2(x,t)+2\theta^2(y)(h(x,t)-h(y,t))^2.
\end{equation}
Since $p_{2,\ep_0}(t)$ is Lipschitz continuous in $(\kappa,\iy)$ and
 $c_0=\ds\sup_{\kappa<t<\iy}|p'_{2,\ep_0}(t)|\le 1/\ep_0$ by (6.5), it follows from (6.6), (6.16), the mean value theorem and Lemma 4.4 that
 \begin{equation}\begin{split}
 J_t(\theta_0,\theta_0)&\le c\,r^{n-2s}+c\iint_{(\BR^n\s B_{5r})\times B_{4r}}\f{1}{|x-y|^{n+2s}}\,dx\,dy \\
 &+ \f{c\,c_0}{r^2}\int_{B_{5r}}\int_{B_{5r}}\f{\int_0^1|\n u(x+\tau(y-x),t)|^2|x-y|^2\,d\tau}{|x-y|^{n+2s}}\,dx\,dy \\
&\le c\bigl(\|u\|^2_{L^{\iy}(\BR^n_I)}+2\bigr)\,r^{n-2s}\le c\,r^{n-2s}.
 \end{split}\end{equation}
By (6.14), (6.15) and (6.17), we have that
\begin{equation}\begin{split}w^{\theta^2_1}_{B_{4r}}(b)\,&\theta^2_1(B_{4r})+\frac{e}{16c}\int_a^b\iint_{B_{2r}\times B_{2r}}\f{|v(x,t)-v(y,t)|^2}{|x-y|^{n+2s}}\,dx\,dy\,dt\\
&\le c\,r^{n}+c\,\kappa^{-1}r^{n} \int_a^b\int_{\BR^n\s B_R}\f{u^-(y,t)}{|y|^{n+2s}}\,dy\,dt+w^{\theta^2_1}_{B_{4r}}(a)\,\theta^2_1(B_{4r}).
\end{split}\end{equation}
For $t\in I$ and $h<d/2$, we denote by
$$\om(t)=|\{x\in B_{r}:u(x,t)\ge d+\kappa\}|\,\,\text{ and }\,\,E^h_t=\{x\in B_{r}:u(x,t)\ge h+\kappa\}.$$
By the assumption (6.2), we have that
$$\int_{-r^{2s}}^0\om(t)\,dt\ge\nu\,|Q_{r}|=\nu\,r^{2s}|B_{r}|.$$
Also, for $\sm\in(0,\nu)$, it is obvious that
$$\int_{-\sm r^{2s}}^0\om(t)\,dt\le\sm r^{2s}|B_{r}|.$$ Thus these two inequalities yield that
\begin{equation}
\int_{-r^{2s}}^{-\sm r^{2s}}\om(t)\,dt\ge(\nu-\sm)\,r^{2s}|B_{r}|.
\end{equation}
By the mean value theorem, there is some $\tau_0\in[-r^{2s},-\sm r^{2s}]$ such that
\begin{equation*}
\om(\tau_0)\ge\f{\nu-\sm}{1-\sm}\,|B_{r}|=\f{\nu-\sm}{2^n(1-\sm)}\,|B_{2r}|.
\end{equation*}
From now on, we take $a=\tau_0$ and $b\in[-\sm r^{2s},0)$ in (6.8).
Since the function $\gm(t):=\ln(2/t)\bigl(\f{(2-t)_+}{(2-t)_++\ep_0}\bigr)^2=r_{2,\ep_0}(t)(2-t)_+^2$ is decreasing in $(0,\iy)$ by Lemma 6.4, we have that
\begin{equation}\begin{split}
w^{\theta^2_1}_{B_{4r}}(b)\,\theta^2_1(B_{4r})&=\int_{B_{4r}}w(x,b)\,\theta^2_1(x)\,dx=\int_{B_{2r}}w(x,b)\,\theta^2_1(x)\,dx \\
&\ge\int_{B_{2r}\s E^h_b}w(x,b)\,\theta^2_1(x)\,dx\ge\theta^2_1(B_{2r}\s E^h_b)\,\gm(2h).
\end{split}\end{equation}
Also we have that
\begin{equation}\begin{split}
w^{\theta^2_1}_{B_{4r}}(\tau_0)\,\theta^2_1(B_{4r})&=\int_{B_{2r}}w(x,\tau_0)\,\theta^2_1(x)\,dx\le\int_{B_{2r}}w(x,\tau_0)\,dx\\
&=\int_{B_{2r}\s E_{\tau_0}^d}w(x,\tau_0)\,dx+\int_{E_{\tau_0}^d}w(x,\tau_0)\,dx \\
&\le\bigl(|B_{2r}|-\om(\tau_0)/2\bigr)\,\gm(h)+\om(\tau_0)\bigl(-\gm(h)/2+\gm(d)\bigr)\\
&\le\biggl(1-\f{\nu-\sm}{2^{n+1}(1-\sm)}\biggr)\gm(h)|B_{2r}|-\,\f{\gm(h)-2\gm(d)}{2}\,\om(\tau_0).
\end{split}\end{equation}
From (6.20) and (6.21), there is a constant $c_1=c_1(n,s,\ld,\Ld)>0$ such that
\begin{equation}\begin{split}\theta^2_1(B_{2r}\s E^h_b)&\le \f{c_1\ds\biggl(1+\kappa^{-1} \biggl(\f{r}{R}\biggr)^{2s}\cT_R(u^-;{\bf 0})\biggr)+\ds\biggl(1-\f{\nu-\sm}{2^{n+1}(1-\sm)}\biggr)\gm(h)}{\gm(2h)}\,|B_{2r}| \\
&\qquad\qquad-\,\f{\gm(h)-2\gm(d)}{2\gm(2h)}\,\om(\tau_0)
\end{split}\end{equation} for any $h\in(0,d/2)$. Here we may choose a sufficiently small $h\in(0,d/2)$ so that $$\theta^2_1(B_{2r}\s E^h_b)\le\vartheta_0\,|B_{2r}|$$ with some constant $\vartheta_0\in(0,1)$.
Thus it follows from (6.18), (6.21) and (6.22) that
\begin{equation}
\int_{\tau_0}^b\iint_{B_{2r}\times B_{2r}}\f{|v(x,t)-v(y,t)|^2}{|x-y|^{n+2s}}\,dx\,dy\,dt\le c_2 \,r^n
\end{equation} for any $b\in[-\sm r^{2s},0)$.

On the other hand, for any $b\in[-\sm r^{2s},0)$ and and $B^0_{5r}\subset B^0_R$, we use the function $$\vp_0(x,t)=\bigl(v_{B_{2r}}^{\theta}(\tau_0)-v_{B_{2r}}^{\theta}(b)\bigr)\f{\theta(x)}{\sqrt\kappa\,u(x,t)}$$ as another testing function.
From (1.3), we obtain that
\begin{equation}\int_{\Om}\pa_t u(x,t)\,\vp_0(x,t)\,dx+\cI_t(u,\vp_0)=0
\end{equation} for $\aee$ $t\in I$.
As in (6.8), the integral of the first term on $[\tau_0,b]$ in (6.24) can be estimated as
\begin{equation}\begin{split}
\int_{\tau_0}^b\int_{\Om}\pa_t u(x,t)\,\vp_0(x,t)\,dx\,dt
&=\f{\bigl(v_{B_{2r}}^{\theta}(\tau_0)-v_{B_{2r}}^{\theta}(b)\bigr)^2}{\sqrt\kappa}\int_{B_{2r}}\theta(x)\,dx \\
&\ge c_3\,r^n\bigl(v_{B_{2r}}^{\theta}(\tau_0)-v_{B_{2r}}^{\theta}(b)\bigr)^2.
\end{split}\end{equation}
Also we have that
\begin{equation*}-\cI_t(u,\vp_0)=\bigl(v_{B_{2r}}^{\theta}(b)-v_{B_{2r}}^{\theta}(\tau_0)\bigr)\,\cI_t(u,\vp_1)\,\,\text{ for $\vp_1(x,t)=\f{\theta(x)}{\sqrt\kappa\,u(x,t)}$}.
\end{equation*}
Thus, by (6.24) and (6.25), applying the mean value theorem and Lemma 4.4 as in (6.17) we have that
\begin{equation}\begin{split}
c_3 r^n\bigl(v_{B_{2r}}^{\theta}(\tau_0)-v_{B_{2r}}^{\theta}(b)\bigr)^2&\le\int_{\tau_0}^b \bigl[-\cI_t(u,\vp_0)\bigr]\,dt \\
&\le\bigl|v_{B_{2r}}^{\theta}(b)-v_{B_{2r}}^{\theta}(\tau_0)\bigr|\int_{\tau_0}^b|\cI_t(u,\vp_1)|\,dt\\
&\le c_4\,r^n\bigl|v_{B_{2r}}^{\theta}(b)-v_{B_{2r}}^{\theta}(\tau_0)\bigr|.
\end{split}\end{equation}
This implies that, for any $\sm\in(0,\nu)$, there is a constant $c_5>0$ depending only on $n,s,\ld,\Ld,\nu$ and $\sm$ such that
\begin{equation}\sup_{b\in[-\sm r^{2s},0)}\,\bigl|v_{B_{2r}}^{\theta}(b)-v_{B_{2r}}^{\theta}(\tau_0)\bigr|\le c_5.
\end{equation}

Finally, by Lemma 6.3, we obtain that
\begin{equation}\begin{split}
&\iint_{Q_r^+}|v(x,t)-v_{Q_r^+}|^2\,dx\,dt\le 4\iint_{Q_r^+}|v(x,t)-v^{\theta}_{Q_{2r}^+}|^2\,dx\,dt\\
&\,\,\le 16\iint_{Q_r^+}|v(x,t)-v^{\theta}_{B_{2r}}(t)|^2\,dx\,dt+16\iint_{Q_r^+}|v^{\theta}_{B_{2r}}(t)-v^{\theta}_{Q_{2r}^+}|^2\,dx\,dt.
\end{split}\end{equation}
We observe that $\kappa\le u\le 1$ in $B_{2r}$ and $\gm(t)=\ln(d/t)$ is Lipschitz continuous in $[\kappa,\iy)$, and so we see that $v(\cdot,t)\in H^s(B_{2r})$ for any $t\in[-r^{2s},0)$.
Since $\theta$ is a nonnegative function with $\theta\equiv\kappa$ on $B_r$, by Lemma 6.2 we  have that
\begin{equation}\begin{split}
&\sqrt\kappa\int_{B_r}|v(x,t)-v^{\theta}_{B_{2r}}(t)|^2\,dx\\
&\qquad\quad\le c\,r^{2s}\iint_{B_{2r}\times B_{2r}}|v(x,t)-v(y,t)|^2\bigl(\theta(x)\wedge\theta(y)\bigr)\,d_K(x,y) \\
&\qquad\quad\le c_6\,\sqrt\kappa\,r^{2s}\iint_{B_{2r}\times B_{2r}}\f{|v(x,t)-v(y,t)|^2}{|x-y|^{n+2s}}\,dx\,dy.
\end{split}\end{equation}
Thus this gives that
\begin{equation}\begin{split}
&\iint_{Q_r^+}|v(x,t)-v^{\theta}_{B_{2r}}(t)|^2\,dx\,dt\\
&\qquad\quad\le c_6\,r^{2s}\int_{-\sm r^{2s}}^0\iint_{B_{2r}\times B_{2r}}\f{|v(x,t)-v(y,t)|^2}{|x-y|^{n+2s}}\,dx\,dy\,dt.
\end{split}\end{equation}
Since we know from simple calculation that
$$v^{\theta}_{Q_{2r}^+}=\f{1}{|Q_{2r}^+|}\iint_{Q_{2r}^+}v^{\theta}_{B_{2r}}(\tau)\,dy\,d\tau,$$ we have that
$$v^{\theta}_{B_{2r}}(t)-v^{\theta}_{Q_{2r}^+}=\f{1}{|Q_{2r}^+|}\iint_{Q_{2r}^+}\bigl(v^{\theta}_{B_{2r}}(t)-v^{\theta}_{B_{2r}}(\tau)\bigr)\,dy\,d\tau.$$ So it follows from Jensen's inequality and (6.27) that
$$\bigl|v^{\theta}_{B_{2r}}(t)-v^{\theta}_{Q_{2r}^+}\bigr|^2\le\f{1}{|Q_{2r}^+|}\iint_{Q_{2r}^+}\bigl|v^{\theta}_{B_{2r}}(t)-v^{\theta}_{B_{2r}}(\tau)\bigr|^2\,dy\,d\tau\le 4 c^2_5.$$
Thus we obtain that
\begin{equation}\iint_{Q_r^+}|v^{\theta}_{B_{2r}}(t)-v^{\theta}_{Q_{2r}^+}|^2\,dx\,dt\le 4 c^2_5\,|Q_r^+|.
\end{equation}
By (6.28), (6.30) and (6.31), we conclude that
\begin{equation*}\begin{split}
&\iint_{Q_r^+}|v(x,t)-v_{Q_r^+}|^2\,dx\,dt\\
&\le 16\,c_6\,r^{2s}\int_{-\sm r^{2s}}^0\iint_{B_{2r}\times B_{2r}}\f{|v(x,t)-v(y,t)|^2}{|x-y|^{n+2s}}\,dx\,dy\,dt+64 \,c^2_5\,|Q_r^+|.
\end{split}\end{equation*}
Therefore we complete the proof. \qed

\section{Nonlocal parabolic Harnack inequality}

\begin{lemma} Let $u\in H^1(I;\rX_g(\Om))$ be a weak solution of the equation ${\bf NP}_{\Om_I}(0,g,g)$ such that $u\ge 0$ in $Q^0_R\subset\Om_I$ where $g\in C(\BR^n_{I_*})\cap L^{\iy}(\BR^n_I)$, and let $m\ge 0$. If there is some $\nu\in(0,1)$ such that
\begin{equation}\bigl|(Q^0_r)^+\cap\{u\ge m\}\bigr|\ge\nu\,|(Q^0_r)^+|
\end{equation} for some $r\in(0,R/5)$, then there is a constant $C_0=C_0(n,s,\ld,\Ld)>0$ such that
\begin{equation}
\biggl|(Q^0_r)^+\cap\biggl\{u\le 2\dt m-\f{1}{2}\biggl(\f{r}{R}\biggr)^{2s}\cT_R(u^-;(x_0,t_0))\biggr\}\biggr|\le\f{C_0}{\nu}\ln^{-1}\biggl(\f{1}{2\dt}\biggr)\,|(Q^0_r)^+|
\end{equation} for any $\dt\in(0,1/4)$.
\end{lemma}

\pf For simplicity in writing, we proceed the proof with $(x_0,t_0)=(0,0):=\bf 0$.
We set $\util=u+\kappa$ where
\begin{equation}\kappa=\f{1}{2}\biggl(\f{r}{R}\biggr)^{2s}\cT_R(u^-;\bf 0).
\end{equation}
Then we see that the function $\bu$ given by
\begin{equation*}\bu(x,t)=\f{\util(x,t)}{m+\kappa}
\end{equation*}
is a weak solution of the equation. Set $v=\ln(1/\bu)$ and take any $\sm\in(0,\nu)$. Then it follows from Theorem 6.1 that
\begin{equation*}\int_a^b\iint_{B_{2r}\times B_{2r}}\f{|v(x,t)-v(y,t)|^2}{|x-y|^{n+2s}}\,dx\,dy\,dt\le C\,r^n
\end{equation*}
for any $a,b\in[-\sm r^{2s},0)$, and also we have the estimate
\begin{equation}\begin{split}
&\iint_{Q_r^+}|v(x,t)-v_{Q_r^+}|^2\,dx\,dt\\
&\qquad\le C\,r^{2s}\int_{-\sm r^{2s}}^0\iint_{B_{2r}\times B_{2r}}|v(x,t)-v(y,t)|^2\,d_K(x,y)+C\,|Q_r^+|\\&\qquad\le C\,|Q_r^+|.
\end{split}\end{equation}
Thus, by Schwarz's inequality and (7.4), we easily obtain that
\begin{equation}
\f{1}{|Q_r^+|}\iint_{Q_r^+}|v(x,t)-v_{Q_r^+}|\,dx\,dt\le C,
\end{equation}
and so this yields that
\begin{equation}
\f{1}{|Q_r^+|}\iint_{Q_r^+}\bigl||v(x,t)|-|v|_{Q_r^+}\bigr|\,dx\,dt\le C
\end{equation}
because $\bigl||v(x,t)|-|v|_{Q_r^+}\bigr|\le |v(x,t)-v_{Q_r^+}|$.
For any $\dt\in(0,1/4)$, we define the function $\bv$ by
$$\bv=(v\vee 0)\wedge\ln\biggl(\f{1}{2\dt}\biggr).$$
Observing the fact that
$$\ap\wedge\bt=\f{\ap+\bt-|\ap-\bt|}{2}\,\,\,\text{ and }\,\,\,\ap\vee\bt=\f{\ap+\bt+|\ap-\bt|}{2}$$ for any $\ap,\bt\in\BR$, we can easily derive from (7.5) and (7.6) that
\begin{equation}
\f{1}{|Q_r^+|}\iint_{Q_r^+}|\bv(x,t)-\bv_{Q_r^+}|\,dx\,dt\le C.
\end{equation}
By the definition of $\bv$, we see that
$$Q_r^+\cap\{\bv=0\}=Q_r^+\cap\{\util\ge m+\kappa\}=Q_r^+\cap\{u\ge m\}.$$ Thus it follows from (7.1) that
\begin{equation}\bigl|Q_r^+\cap\{\bv=0\}\bigr|\ge\nu\,|Q_r^+|.
\end{equation}
From (7.8), we can derive the following estimate
\begin{equation}\begin{split}
\ln\biggl(\f{1}{2\dt}\biggr)&=\f{1}{|Q_r^+\cap\{\bv=0\}|}\iint_{Q_r^+\cap\{\bv=0\}}\biggl[\ln\biggl(\f{1}{2\dt}\biggr)-\bv(x,t)\biggr]\,dx\,dt\\
&\le\f{1}{\nu}\,\biggl[\ln\biggl(\f{1}{2\dt}\biggr)-\bv_{Q_r^+}\biggr].
\end{split}\end{equation}
Integrating the inequality (7.9) on $Q_r^+\cap\{\bv=\ln(1/2\dt)\}$ and applying (7.1), we easily obtain that
\begin{equation*}\begin{split}
&\biggl|Q_r^+\cap\biggl\{\bv=\ln\biggl(\f{1}{2\dt}\biggr)\biggr\}\biggr|\,\ln\biggl(\f{1}{2\dt}\biggr) \\
&\qquad\qquad\le\f{1}{\nu}\iint_{Q_r^+}|\bv(x,t)-\bv_{Q_r^+}|\,dx\,dt
\le\f{C_0}{\nu}\,|Q_r^+|.
\end{split}\end{equation*}
This implies that
\begin{equation*}
\bigl|Q_r^+\cap\{\util\le 2\dt m\}\bigr|\le\f{C_0}{\nu}\,\ln^{-1}\biggl(\f{1}{2\dt}\biggr)\,|Q_r^+|.
\end{equation*}
Hence we complete the proof. \qed

\begin{lemma} Let $u\in H^1(I;\rX_g(\Om))$ be a weak solution of the equation ${\bf NP}_{\Om_I}(0,g,g)$ such that $u\ge 0$ in $Q^0_R\subset\Om_I$ where $g\in C(\BR^n_{I_*})\cap L^{\iy}(\BR^n_I)$, and let $m\ge 0$. If there is some $\nu\in(0,1]$ such that
\begin{equation*}\bigl|(Q^0_r)^+\cap\{u\ge m\}\bigr|\ge\nu\,|(Q^0_r)^+|
\end{equation*} for some $r\in(0,R/5)$, then there is a constant $\dt=\dt(n,s,\ld,\Ld)\in(0,1/4)$ such that
\begin{equation*}
\inf_{(Q^0_r)^+}u\ge\dt m-\biggl(\f{r}{R}\biggr)^{2s}\cT_R(u^-;(x_0,t_0)).
\end{equation*}
\end{lemma}

\pf For simplicity, we proceed the proof by setting $(x_0,t_0)=(0,0):=\bf 0$.
Since $u\ge 0$ in $Q_R$, without loss of generality we may assume that
\begin{equation}\dt m\ge\biggl(\f{r}{R}\biggr)^{2s}\cT_R(u^-;\bf 0).
\end{equation}
Set $w=(h-u)_+$ for $h\in(\dt m,2\dt m)$. For $\vr\in(r,2r)$, we now choose a testing function $\vp(x,t)=\theta^2(x)\e^2(t) w(x,t)$ where $\theta\in C_c^{\iy}(B_{\vr})$ is a function satisfying that $0\le\theta\le 1$ and $|\n\theta|\le c/\vr$ in $\BR^n$, and $\e\in C^{\iy}_c(-\sm\vr^{2s},\iy]$ is a function such that $0\le\e\le 1$ and $|\e'|\le c/\vr^{2s}$ in $\BR$. Then we have that
\begin{equation}\begin{split}
0&=\iint_{\BR^n\times\BR^n}(u(x,t)-u(y,t))\bigl(\vp(x,t)-\vp(y,t)\bigr)K(x-y)\,dx\,dy\\
&\qquad\qquad\qquad\qquad\qquad+\int_{\Om}\pa_t u(x,t)\,\vp(x,t)\,dx\\
&=\int_{B_{\vr}}\int_{B_{\vr}}(u(x,t)-u(y,t))\bigl(\vp(x,t)-\vp(y,t)\bigr)K(x-y)\,dx\,dy\\
&\qquad+2\int_{\BR^n\s B_{\vr}}\int_{B_{\vr}}(u(x,t)-u(y,t))\,\vp(x,t)K(x-y)\,dx\,dy\\
&\qquad\qquad\qquad\qquad\qquad+\int_{B_{\vr}}\pa_t u(x,t)\,\vp(x,t)\,dx\\
&:=\cI_1(t)+\cI_2(t)+\cI_3(t)
\end{split}\end{equation} for $\aee$ $t\in I$. Splitting $\cI_2(t)$ into two parts yields that
\begin{equation*}\begin{split}
\f{1}{2}\,\cI_2(t)&=\int_{B^c_{\vr}\cap\{y:u(y,t)<0\}}\int_{B_{\vr}}(u(x,t)-u(y,t))\,\vp(x,t)K(x-y)\,dx\,dy\\
&+\int_{B^c_{\vr}\cap\{y:u(y,t)\ge 0\}}\int_{B_{\vr}}(u(x,t)-u(y,t))\,\vp(x,t)K(x-y)\,dx\,dy\\
&=\cI_2^1(t)+\cI_2^2(t).
\end{split}\end{equation*}
Thus we have that
\begin{equation*}\begin{split}
\cI_2^1(t)&\le h|B_{\vr}\cap\{x:u(x,t)<h\}|\sup_{x\in\supp(\theta)}\int_{B^c_{\vr}}(h+u^-(y,t))K(x-y)\,dy,\\
\cI_2^2(t)&\le h^2|B_{\vr}\cap\{x:u(x,t)<h\}|\sup_{x\in\supp(\theta)}\int_{B^c_{\vr}}K(x-y)\,dy,
\end{split}\end{equation*}
and this leads us to get that
\begin{equation}\begin{split}&\int^0_{-\sm\vr^{2s}}\cI_2(t)\,dt\\
&\quad\le 4h\,|Q_{\vr}^+\cap\{u<h\}|
\sup_{(x,t)\in\Gm^{\theta}_{\vr}}\int_{B^c_{\vr}}(h+u^-(y,t))K(x-y)\,dy,
\end{split}\end{equation}
where $\Gm^{\theta}_{\vr}=\supp(\theta)\times[-\sm\vr^{2s},0)$.
As in the proof of the Caccioppoli type estimate in Theorem 1.4 \cite{DKP1}, we have that
\begin{equation}\begin{split}
&\int^0_{-\sm\vr^{2s}}\cI_1(t)\,dt\\
&\le-c\int^0_{-\sm\vr^{2s}}\iint_{B_{\vr}\times B_{\vr}}|w(x,t)\theta(x)-w(y,t)\theta(y)|^2 K(x-y)\,dx\,dy\,dt\\
&+c\int^0_{-\sm\vr^{2s}}\iint_{B_{\vr}\times B_{\vr}}(w(x,t)\vee w(y,t))^2|\theta(x)-\theta(y)|^2 K(x-y)\,dx\,dy\,dt.
\end{split}\end{equation}
For any $\tau\in[-\sm\vr^{2s},0)$, we have the estimate
\begin{equation}\begin{split}
\int_{-\sm\vr^{2s}}^{\tau}\cI_3(t)\,dt&=-\int_{B_{\vr}}w^2(x,\tau)\theta^2(x)\e^2(\tau)\,dx\\
&\quad+2\int_{-\sm\vr^{2s}}^\tau\int_{B_{\vr}} w^2(x,t)\theta^2(x)\e(t)\e'(t)\,dx\,dt.
\end{split}\end{equation}
From (7.11), (7.12), (7.13) and (7.14), we have that
\begin{equation}\begin{split}
&\sup_{t\in[-\sm\vr^{2s},0)}\e^2(t)\|w(\cdot,t)\theta\|^2_{L^2(B_{\vr})}
+c\int_{-\sm\vr^{2s}}^0\e^2(t)[w(\cdot,t)\theta]^2_{H^s(B_{\vr})}dt\\
&\,\le 4h\,|Q_{\vr}^+\cap\{u<h\}|
\sup_{(x,t)\in\Gm^{\theta}_{\vr}}\int_{B^c_{\vr}}(h+u^-(y,t))K(x-y)\,dy\\
&\,+c\int^0_{-\sm\vr^{2s}}\iint_{B_{\vr}\times B_{\vr}}(w(x,t)\vee w(y,t))^2|\theta(x)-\theta(y)|^2\,K(x-y)\,dx\,dy\,dt\\
&\,+c\vr^{-2s}\|w\|^2_{L^2(Q_{\vr}^+)}:= 4h\,|Q_{\vr}^+\cap\{u<h\}|\,A_{\vr}(u,\theta,h)+B_{\vr}(w,\theta)+C_{\vr}(w).
\end{split}\end{equation}
We now apply Lemma 2.2 to the estimates (7.15). For $k=0,1\cdots$, we set
\begin{equation}h_k=\dt m+2^{-(k+1)}\dt m,\,\,\,\,\vr_k=r+2^{-k}r,\,\,\,\,\overline\vr_k=\f{\vr_{k+1}+\vr_k}{2}.
\end{equation}
We observe that $r<\vr_k,\overline\vr_k<2r$ and $h_k-h_{k+1}\ge 2^{-(k+3)}h_k$, and
\begin{equation*}h_0=\f{3}{2}\dt m\le 2\dt m-\f{1}{2}\biggl(\f{r}{R}\biggr)^{2s}\cT_R(u^-;\bf 0).
\end{equation*}
Since we see that
$$\{u<h_0\}\subset\biggl\{u\le 2\dt m-\f{1}{2}\biggl(\f{r}{R}\biggr)^{2s}\cT_R(u^-;\bf 0)\biggr\},$$
by Lemma 7.1 we have that
\begin{equation}\f{|Q_r^+\cap\{u<h_0\}|}{|Q_r^+|}\le\f{C_0}{\nu}\ln^{-1}\biggl(\f{1}{2\dt}\biggr).
\end{equation}
For $k=0,1,\cdots$, we note that
\begin{equation}\begin{split}w_k&=(h_k-u)_+ \ge(h_k-h_{k+1})\mathbbm{1}_{\{u<h_{k+1}\}} \\
&\ge 2^{-k-3}h_k\mathbbm{1}_{\{u<h_{k+1}\}}.
\end{split}\end{equation}
For $k=0,1,\cdots$, set $Q_k=Q_{\vr_k}^+$, and let $\theta_k\in C_c^{\iy}(B_{\overline\vr_k})$ be a function with $\theta_k\equiv 1$ on $B_{\vr_{k+1}}$ such that $0\le\theta_k\le 1$ and $|\n\theta_k|\le 2^{k+1}/r$ in $\BR^n$ and let $\e_k\in C^{\iy}_c(-\sm\overline\vr^{2s}_k,\iy]$ is a function with $\e_k\equiv 1$ on $[-\sm\vr_{k+1}^{2s},0)$ such that $0\le\e_k\le 1$ and
$$|\e'_k|\le\f{2^{k+1}}{\sm\overline\vr_k^{2s}}\,\,\text{ in $\BR$.}$$
By applying H\'older's inequality with $$q=\f{n}{n-2s}\,\,\text{ and }\,\,q'=\f{n}{2s},$$ and the fractional Sobolev's inequality with $$\ap=1+\f{2s}{n}=\f{1}{q'}+1,$$ we can derive the following inequalities
\begin{equation}\begin{split}
&\iint_{Q_{k+1}}|w_k|^{2\ap}\,dx\,dt\le\iint_{Q_k}|w_k\theta_k\e_k|^2|w_k\theta_k\e_k|^{\f{4s}{n}}\,dx\,dt\\
&\quad\le\int_{-\sm\vr_k^{2s}}^0\biggl(\int_{B_{\vr_k}}|w_k\theta_k\e_k|^{2q}\,dx\biggr)^{\f{1}{q}}\biggl(\int_{B_{\vr_k}}|w_k\theta_k\e_k|^2\,dx\biggr)^{\f{1}{q'}}dt\\
&\quad\le\gm_{n,2}\biggl(\sup_{\,t\in[-\sm\vr^{2s}_k,0)}\e_k^2(t)\|w_k(\cdot,t)\theta_k\|^2_{L^2(B_{\vr_k})}\biggr)^{\f{1}{q'}}\\
&\qquad\times\biggl[(1-s)\int^0_{-\sm\vr^{2s}_k}\e_k^2(t)[w_k(\cdot,t)\theta_k]^2_{H^s(B_{\vr_k})}\,dt+\vr_k^{-2s}\|w_k\|^2_{L^2(Q_k)}\biggr]\\
&\quad\le C\biggl(4h_k\,|Q_k\cap\{u<h_k\}|\,A_{\vr_k}(u,h_k)+B_{\vr_k}(w_k,\theta_k)+C_{\vr_k}(w_k)\biggr)^{\ap}
\end{split}\end{equation}
where $C>0$ is a universal constant depending only on $n$ and $s$. From (7.10), (7.15) and the fact that
\begin{equation*}
|y-x|\ge|y|-|x|\ge\biggl(1-\f{\bar\vr_k}{\vr_k}\biggr)|y|\ge 2^{-k-3}|y|
\end{equation*}
for any $y\in B^c_{\vr_k}$ and $x\in B_{\bar\vr_k}$, by (7.10) and (7.16) we get the estimate
\begin{equation}\begin{split}
A_{\vr_k}(u,\theta_k,h_k)&=\sup_{(x,t)\in\Gm^{\theta_k}_{\vr_k}}\int_{B^c_{\vr_k}}(h_k+u^-(y,t))K(x-y)\,dy \\
&\le c\,2^{k(n+2s)}\int_{B^c_{\vr_k}}\f{h_k+u^-(y,t)}{|y|^{n+2s}}\,dy \\
&\le c\,2^{k(n+2s)}h_k\,r^{-2s}+c\,2^{k(n+2s)}\int_{B^c_R}\f{u^-(y,t)}{|y|^{n+2s}}\,dy\\
&\le c\,2^{k(n+2s)}h_k\,r^{-2s}+c\,2^{k(n+2s)} r^{-2s}\biggl(\f{r}{R}\biggr)^{2s}\cT_R(u^-;\bf 0) \\
&\le c\,2^{k(n+2s)}h_k\,r^{-2s}.
\end{split}\end{equation}
Also we have the following estimates
\begin{equation}\begin{split}
&B_{\vr_k}(w_k,\theta_k)\\&\le c\int^0_{-\sm\vr_k^{2s}}\iint_{B_{\vr_k}\times B_{\vr_k}}(w_k(x,t)\vee w_k(y,t))^2\,\f{|\theta_k(x)-\theta_k(y)|^2}{|x-y|^{n+2s}}\,dx\,dy\,dt\\
&\le c\,h^2_k\int^0_{-\sm\vr_k^{2s}}\int_{B_{\vr_k}}\int_{B_{\vr_k}\cap\{u(\cdot,t)<h_k\}}\f{\sup_{\BR^n}|\n\theta_k|^2}{|x-y|^{n+2s-2}}\,dx\,dy\,dt\\
&\le c\,h_k^2\,\biggl(\f{2^k}{r}\biggr)^2\int^0_{-\sm\vr_k^{2s}}\int_{B_{\vr_k}\cap\{u(\cdot,t)<h_k\}}\biggl(\,\int_{B_{2\vr_k}}\f{1}{|y|^{n+2s-2}}\,dy\biggr)dx\,dt\\
&\le c\,2^{2k}\,h_k^2\,r^{-2s}\,|Q_k\cap\{u<h_k\}|.
\end{split}\end{equation}
and $\,C_{\vr_k}(w_k)\le c\,h_k^2\,r^{-2s}\,|Q_k\cap\{u<h_k\}|.$ From (7.19), (7.20) and (7.21), we conclude that
$$\biggl(\iint_{Q_{k+1}}|w_k|^{2\ap}\,dx\,dt\biggr)^{\f{1}{\ap}}\le C\,2^{k(n+2s+2)}\,h_k^2\,r^{-2s}\,|Q_k\cap\{u<h_k\}|.$$
Since $|Q_{k+1}|^{-1/\ap}\sim r^{-n}$ and $|Q_k|\sim r^{n+2s}$, this estimate and (7.18) yield that
\begin{equation}\begin{split}
&(h_k-h_{k+1})^2\biggl(\f{|Q_{k+1}\cap\{u<h_{k+1}\}|}{|Q_{k+1}|}\biggr)^{\f{1}{\ap}}\\&\qquad\qquad\qquad\qquad\qquad\le\biggl(\f{1}{|Q_{k+1}|}\iint_{Q_{k+1}}|w_k|^{2\ap}\,dx\,dt\biggr)^{\f{1}{\ap}}\\
&\qquad\qquad\qquad\qquad\qquad\le C\,2^{k(n+2s+2)}\,h_k^2\,\f{|Q_k\cap\{u<h_k\}|}{|Q_k|}.
\end{split}\end{equation}
For $k=0,1,2,\cdots,$ we set
$$N_k=\f{|Q_k\cap\{u<h_k\}|}{|Q_k|}.$$
Then it follows from (7.22) that
$$N_{k+1}^{\f{1}{\ap}}\le C\,\f{2^{k(n+2s+2)}\,h_k^2}{(h_k-h_{k+1})^2}\,N_k
\le C\,2^{k(n+2s+4)}N_k.$$ This leads to us to obtain that
\begin{equation}
N_{k+1}\le C_1\,2^{k\ap(n+2s+4)}N_k^{1+\f{2s}{n}}
\end{equation} where $C_1=C^{\ap}$.
In addition, we see from (7.17) that
\begin{equation}
N_0\le\f{C_0}{\nu}\,\ln^{-1}\biggl(\f{1}{2\dt}\biggr).
\end{equation}
We apply Lemma 2.2 with
$$d_0=C_1,\,\,\,a=2^{n+2s+4}>1,\,\,\,\text{ and }\,\,\,\e=\f{2s}{n}.$$
If we choose a small $\dt$ depending only on $n,s,\ld,\Ld$ and $\nu$ so that
\begin{equation}0<\dt:=\f{1}{2}\exp\biggl(-\f{C_0 C_1^{\f{n}{2s}}2^{\f{n^2}{4 s^2}(n+2s+4)}}{\nu}\biggr)<\f{1}{4},
\end{equation}
then $N_0\le C_1^{-\f{n}{2s}}(2^{n+2s+4})^{-\f{n^2}{4 s^2}}$. Thus we conclude that $\lim_{k\to\iy}N_k=0$. This implies that
$$\inf_{Q_r^+}u\ge\dt m.$$ Hence we comlete the proof.
\qed

\,\,Next, we need a parabolic version of the Krylov-Safonov covering theorem \cite{KS} which is a useful tool for the proof of Theorem 7.4.

For $(x,t)\in\BR^n\times\BR$, $r>0$ and $\sm\in(0,\nu)$ (where $\nu$ is the constant given in (6.3)), we see that the parabolic cylinders
$$Q^+_r(x,t)=B_r(x)\times[t-\sm r^{2s},t)$$ is given by $Q^+_r(x,t)=\{Y\in\BR^n\times\BR:\dd(Y,X)<r\}$ where $\dd$ is
the parabolic distance between $X=(x,t)$ and $Y=(y,\tau)$ by
\begin{equation}\dd(X,Y)=\begin{cases} |x-y|\vee\biggl(\ds\f{|t-\tau|}{\sm}\biggr)^{1/2s}, & \tau<t,\\
                          \iy, & \tau\ge t. \end{cases}
\end{equation}
Then we note that  If $E\subset\BR^n\times\BR$ is a bounded set and $\fC_E$ is a collection of cylinders $Q^+_r(x,t)$ with $(x,t)\in E$, then it follows from Vitali's covering theorem that there is a countable pairwise disjoint subcollection $\cC_E=\{Q^+_{r_k}(x_k,t_k)\}_{k\in\BN}$ of $\fC_E$ such that
$$E\subset\bigcup_{k\in\BN}Q^+_{r_k}(x_k,t_k).$$
For $\vr>0$, $\gm\in(0,1)$ and a measurable subset $E$ of a cylinder $(Q^0_r)^+=Q^+_r(x_0,t_0)$, we define the set
\begin{equation}E^{\vr}_{\gm}=\bigcup_{0<\rho<\vr}\bigl\{Q^+_{3\rho}(X)\cap(Q^0_r)^+:|E\cap Q^+_{3\rho}(X)|>\gm\,|Q^+_{\rho}(X)|,X\in(Q^0_r)^+\bigr\}.
\end{equation}
The following nonlocal parabolic version of the Krylov-Safonov covering theorem no longer depends on the threshold radius $\vr$. Its proof is based on that of \cite{KSh}.

\begin{lemma} If $\vr>0$, $\gm\in(0,1)$ and $E\subset Q^+_r(x_0,t_0)$ is a measurable set, then

\,either $\,\,|E^{\vr}_{\gm}|\ds\ge\f{2^{-(n+2s)}}{\gm}\,|E|\,\,$ or $\,\,E^{\vr}_{\gm}=Q^+_r(x_0,t_0).$
\end{lemma}

\,\pf We define the maximal operator $\cM:\BR^n\times\BR\to\BR$ by
$$\cM(y,\tau)=\sup_{Q^+_{3\rho}(x,t)\in\fF}\f{|E\cap Q^+_{3\rho}(x,t)|}{|Q^+_{\rho}(x,t)|}$$
where $\fF$ is the family of all cylinders $Q^+_{3\rho}(x,t)$ with $(x,t)\in Q^+_r(x_0,t_0)$, $0<\rho<\vr$ and $(y,\tau)\in Q^+_{3\rho}(x,t)$. Then we see that
\begin{equation}E^{\vr}_{\gm}=\{(y,\tau)\in Q^+_r(x_0,t_0):\cM(y,\tau)>\gm\}.
\end{equation}
Indeed, if $Y=(y,\tau)\in Q^+_r(x_0,t_0)$, then there is a cylinder $Q^+_{3\rho}(x,t)$
with $(x,t)\in Q^+_r(x_0,t_0)$, $0<\rho<\vr$ and $|E\cap Q^+_{3\rho}(x,t)|>\gm|Q^+_{\rho}(x,t)|$. Thus this gives that $Y\in E^{\vr}_{\gm}$. On the other hand, if $Y\in E^{\vr}_{\gm}$, then there is a cylinder $Q^+_{3\rho}(x,t)$
with $(x,t)\in Q^+_r(x_0,t_0)$ and $Y\in Q^+_{3\rho}(x,t)$, $0<\rho<\vr$ and $|E\cap Q^+_{3\rho}(x,t)|>\gm|Q^+_{\rho}(x,t)|$. This implies that $\cM(Y)>\gm$.

Suppose that $E^{\vr}_{\gm}\neq Q^+_r(x_0,t_0)$, i.e. $Q^+_r(x_0,t_0)\s E^{\vr}_{\gm}\neq\phi$. Since $E^{\vr}_{\gm}$ is open with respect to the metric $\dd$,
we see that
$$\rho_Y:=\f{1}{2}\sup\{\rho>0:Y\in Q^+_{\rho}(x,t)\subset E^{\vr}_{\gm}, Q^+_{2\rho}(x,t)\cap(Q^+_r(x_0,t_0)\s E^{\vr}_{\gm})=\phi\}>0$$
for each $Y=(y,\tau)\in E^{\vr}_{\gm}$. So we may assume that $\rho_Y<\vr/4$. Since each $Y\in E^{\vr}_{\gm}$ yields a point $X_Y:=(x_Y,t_Y)$ so that $Y\in Q^+_{\rho_Y}(X_Y)\subset E^{\vr}_{\gm}$ and
$$Q^+_{5\rho_Y}(X_Y)\cap(Q^+_r(x_0,t_0)\s E^{\vr}_{\gm})\neq\phi,$$ the family $\cC=\{Q^+_{\rho_Y}(X_Y):Y\in E^{\vr}_{\gm}\}$ covers $E^{\vr}_{\gm}$. From Vitali's covering lemma, we can extract a countable family $\{Q^+_{\rho_k}(X_k)\}_{k\in\BN}$ of pairwise disjoint parabolic cubes (where $X_k=X_{Y_k}$ and $\rho_k=\rho_{Y_k}$ for $k\in\BN$) such that
$$E^{\vr}_{\gm}\subset\bigcup_{k\in\BN}Q^+_{\rho_k}(X_k).$$
Observe that $Q^+_{5\rho_k}(X_k)\cap(Q^+_r(x_0,t_0)\s E^{\vr}_{\gm})\neq\phi$
for all $k\in\BN$. If $Y_k=(y_k,\tau_k)\in Q^+_{5\rho_k}(X_k)\cap(Q^+_r(x_0,t_0)\s E^{\vr}_{\gm})$, then we easily see that $\cM(Y_k)\le\gm$ for $k\in\BN$.
Since $Y_k\in Q^+_{5\rho_k}(X_k)$ and $5\rho_k/3<\vr$, we have that
\begin{equation}|E\cap Q^+_{5\rho_k}(X_k)|\le\gm|Q^+_{5\rho_k/3}(X_k)|\le 2^{n+2s}\gm|Q^+_{\rho_k}(X_k)|.
\end{equation}
Moreover, by (7.28) we see that every density points of $E$ belongs to $E^{\vr}_{\gm}$, because
$$\liminf_{\rho\to 0}\f{|E\cap Q^+_{3\rho}(X)|}{|Q^+_{\rho}(X)|}\ge 1>\gm$$
for any density point $X$ of $E$. Hence this and (7.29) enables us to obtain that
$$|E|=|E\cap E^{\vr}_{\gm}|\le\sum_{k=1}^{\iy}|E\cap Q^+_{5\rho_k}(X_k)|\le 2^{n+2s}\gm\sum_{k=1}^{\iy}|Q^+_{\rho_k}(X_k)|\le 2^{n+2s}\gm|E^{\vr}_{\gm}|.\qed$$

\begin{thm} Let $g\in C(\BR^n_{I_*})\cap L^{\iy}(\BR^n_I)$. If $u\in H^1(I;\rX_g(\Om))$ is a weak solution of the nonlocal equation ${\bf NP}_{\Om_I}(0,g,g)$ with $u\ge 0$ in $Q^0_R\subset\Om_I$, then we have the estimate
\begin{equation}\biggl(\f{1}{2|(Q^0_r)^+|}\int_{(Q^0_r)^+}u^p\,dx\,dt\biggr)^{\f{1}{p}}\le \inf_{(Q^0_r)^+}u+\f{4}{3}\,\biggl(\f{r}{R}\biggr)^{2s}\cT_r(u^-;(x_0,t_0))
\end{equation} for any $p\in(0,1)$ and $r\in (0,R)$.
\end{thm}

\pf Take any $r\in(0,R)$.  For simplicity, we may assume that $(x_0,t_0)=(0,0)$.
  For $\ap>0$ and $k=0,1,2,\cdots,$ we set
$$\cR^k_\ap=\biggl\{(x,t)\in Q^+_r:u(x,t)\ge\dt^k\ap-\f{2\kappa}{1-\dt}\biggr\}$$
where $\dt\in(0,1/4)$ is the constant in Lemma 7.2 and $\kappa>0$ is the constant given by
$$\kappa=\f{1}{2}\biggl(\f{r}{R}\biggr)^{2s}\cT_R(u^-;(0,0)).$$
Then we see that $\cR^{k-1}_{\ap}\subset\cR^k_{\ap}$ for all $k\in\BN$.
Let $(x,t)\in Q^+_r$ with $Q^+_{3\rho}(x,t)\cap Q^+_r\subset E^{\vr}_{\gm}$ where $E=\cR^{k-1}_{\ap}$. From (7.27), we see that
$$|\cR^{k-1}_{\ap}\cap Q^+_{3\rho}(x,t)|>\gm|Q^+_{\rho}(x,t)|=\f{\gm}{3^{n+2s}}|Q^+_{3\rho}(x,t)|.$$
Applying Lemma 7.2 with $m=\dt^{k-1}\ap-\f{2\kappa}{1-\dt}$ and $\nu=\gm\,3^{-n-2s}$, we have that
$$u\ge\dt\biggl(\dt^{k-1}\ap-\f{2\kappa}{1-\dt}\biggr)-2\kappa=\dt^k\ap-\f{2\kappa}{1-\dt}\,\,\,\text{ in $Q^+_{3\rho}(x,t)\cap Q^+_r$}$$
and so we obtain that $E^{\vr}_{\gm}\subset\cR^k_{\ap}$. Thus it follows from Lemma 7.3 that either $\cR^k_{\ap}=Q^+_r$ or $|\cR^k_{\ap}|\ge\f{2^{-n-2s}}{\gm}|\cR^{k-1}_{\ap}|$ for any $k\in\BN$. Here, without loss of generality, we assume that
\begin{equation}\f{1}{16}\,2^{-n-2s}<\gm<2^{-n-2s}.
\end{equation}
Then we claim that, if there is some $N\in\BN$ such that
\begin{equation}|\cR^0_{\ap}|>(2^{n+2s}\gm)^N|Q^+_r|,
\end{equation}
then we have that $\cR^N_{\ap}=Q^+_r$. Indeed, if $\cR^N_{\ap}\neq Q^+_r$, then $|\cR^N_{\ap}|\ge\f{2^{-n-2s}}{\gm}|\cR^{N-1}_{\ap}|$, and so
$\cR^k_{\ap}\neq Q^+_r$ for any $k=1,2,\cdots,N.$ This implies that
$$|\cR^N_{\ap}|\ge\f{1}{2^{n+2s}\gm}\,|\cR^{N-1}_{\ap}|\ge\cdots\ge\biggl(\f{1}{2^{n+2s}\gm}\biggr)^N|\cR^0_{\ap}|>|Q^+_r|,$$
which gives a contradiction. Thus the fact that $\cR^N_{\ap}=Q^+_r$ leads us to obtain that
\begin{equation}u\ge\dt^N\ap-\f{2\kappa}{1-\dt}\,\,\,\text{ in $Q^+_r$.}
\end{equation}
If $N$ is the smallest integer satisfying (7.31), we see that
\begin{equation}N>\f{1}{\ln(2^{n+2s}\gm)}\ln\biggl(\f{|\cR^0_{\ap}|}{|Q^+_r|}\biggr).
\end{equation}
From (7.32) and (7.33), we conclude that
\begin{equation}\inf_{Q^+_r}u\ge\ap\biggl(\f{|\cR^0_{\ap}|}{|Q^+_r|}\biggr)^{1/\el}-\f{2\kappa}{1-\dt}\,\,\,\text{ with $\el=\ds\f{\ln(2^{n+2s}\gm)}{\ln\dt}$,}
\end{equation}
where $\dt$ and $\el$ depend only on $n,s,\ld$ and $\Ld$. This enables us to get that
\begin{equation}\begin{split}
\f{|Q^+_r\cap\{u\ge\ap-\f{2\kappa}{1-\dt}\}|}{|Q^+_r|}=\f{|\cR^0_{\ap}|}{|Q^+_r|} \le\ap^{-\el}\biggl(\,\inf_{Q^+_r}u+\f{2\kappa}{1-\dt}\biggr)^{\el}.
\end{split}\end{equation}
By standard analysis, we have that
\begin{equation}\f{1}{|Q^+_r|}\iint_{Q^+_r}u^p\,dx\,dt=p\int_0^{\iy}\ap^{p-1}\f{|Q^+_r\cap\{u\ge\ap\}|}{|Q^+_r|}\,d\ap
\end{equation} for any $p>0$. Thus it follows from (7.35) and (7.36) that
\begin{equation*}\begin{split}
\f{1}{|Q^+_r|}\iint_{Q^+_{2r}}u^p\,dx\,dt&\le p\int_0^{\iy}\ap^{p-1}\f{|Q^+_r\cap\{u\ge\ap-\f{2\kappa}{1-\dt}\}|}{|Q^+_r|}\,d\ap\\
&\le p\int_0^h\ap^{p-1}\,d\ap+p\biggl(\,\inf_{Q^+_r}u
+\f{2\kappa}{1-\dt}\biggr)^{\el}\int_h^{\iy}\ap^{p-1-\el}\,d\ap.
\end{split}\end{equation*}
If we take $h=\ds\inf_{Q^+_r}u
+\f{2\kappa}{1-\dt}$ and $p=\el/2$, then we conclude that
$$\f{1}{|Q^+_r|}\iint_{Q^+_r}u^p\,dx\,dt\le 2\biggl(\,\inf_{Q^+_r}u
+\f{2\kappa}{1-\dt}\biggr)^p\le 2\biggl(\,\inf_{Q^+_r}u
+\f{4}{3}\,2\kappa\biggr)^p$$
and $0<p<1$ by (7.30) because $\dt\in(0,1/4)$.
Since we can take any sufficiently large $C_0, C_1$ in (7.25), by (7.30) and (7.34) the result of Theorem 7.4
holds for all $p\in(0,1)$ with the same universal constants $1$ and $4/3$.
Hence we are done. \qed

\begin{thm} Let $g\in C(\BR^n_{I_*})\cap L^{\iy}(\BR^n_I)$. If $u\in H^1(I;\rX_g(\Om))$ is a weak solution of the nonlocal parabolic equation ${\bf NP}_{\Om_I}(0,g,g)$ with $u\ge 0$ in $Q^0_R\subset\Om_I$, then there exists a constnat $c>0$ depending only on $n,s,\ld$ and $\Ld$ such that
\begin{equation}\sup_{(Q^0_r)^-}u\le c\,\inf_{(Q^0_r)^+}u+c\,\biggl(\f{r}{R}\biggr)^{2s}\cT_r(u^-;(x_0,t_0))
\end{equation} for any $r>0$ with $5r<R$.
\end{thm}

We can easily derive the following nonlocal parabolic Harnack inequalities for a nonnegative weak solutions of the nonlocal parabolic equation ${\bf NP}_{\Om_I}(0,g,g)$ as a natural by-product of Theorem 7.5. By the way, the result has no nonlocal parabolic tail, which means that the result coincides with that of local parabolic case.

\begin{cor} Let $g\in C(\BR^n_{I_*})\cap L^{\iy}(\BR^n_I)$. If $u\in H^1(I;\rX_g(\Om))$ is any nonnegative weak solution of the nonlocal parabolic equation ${\bf NP}_{\Om_I}(0,g,g)$, then there exists a constnat $c>0$ depending only on $n,s,\ld$ and $\Ld$ such that
\begin{equation}\sup_{(Q^0_r)^-}u\le c\,\inf_{(Q^0_r)^+}u
\end{equation} for any $r>0$ with $5r<R$.
\end{cor}

\,{\bf [Proof of Theorem 7.5(nonlocal Harnack inequality)]}   If $Q^0_{3 r}\subset Q^0_R\subset\Om_I$, then by Theorem 5.2 and Lemma 5.3 there is a constant $c_0=c_0(n,s,\ld,\Ld)>0$ such that
\begin{equation}\begin{split}
\sup_{Q^0_{r/2}}u&\le\dt\,\cT_r(u^+;(x_0,t_0))+c_0\,\dt^{-\gm_s}\biggl(\displaystyle\f{1}{|\cQ^0_{2r}|}\iint_{Q^0_r}u^2\,dx\,dt\biggr)^{\f{1}{2}}  \\
&\le c_0\,\dt^{-\gm_s}\biggl(\f{1}{|\cQ^0_{2r}|}\iint_{Q^0_{2r}}u^2\,dx\,dt\biggr)^{\f{1}{2}} \\
&\qquad\qquad +c\,\dt\sup_{Q^0_{2r}}u+c\,\dt\biggl(\f{r}{R}\biggr)^{2s}\cT_R(u^-;(x_0,t_0))
\end{split}\end{equation}
for any $\dt\in(0,1]$, where $\gm_s=\f{n+2s}{4s}$.
Thus it follows from a covering argument with $\f{1}{2}\le a<b\le 2$ that
\begin{equation}\begin{split}
\sup_{Q^0_{ar}}u&\le \f{c\,\dt^{-\gm_s}}{(b-a)^{\f{n+2s}{2}}}\biggl(\f{1}{|\cQ^0_{br}|}\iint_{Q^0_{br}}u^2\,dx\,dt\biggr)^{\f{1}{2}} \\
&\qquad\qquad +c\,\dt\sup_{Q^0_{br}}u+c\,\dt\biggl(\f{r}{R}\biggr)^{2s}\cT_R(u^-;(x_0,t_0)) \\
&\le\f{c\,\dt^{-\gm_s}}{(b-a)^{\f{n+2s}{2}}}\biggl(\sup_{Q^0_{br}}u\biggr)^{\f{2-p}{2}}\biggl(\f{1}{|\cQ^0_{br}|}\iint_{Q^0_{br}}u^p\,dx\,dt\biggr)^{\f{1}{2}} \\
&\qquad\qquad +c\,\dt\sup_{Q^0_{br}}u+c\,\dt\biggl(\f{r}{R}\biggr)^{2s}\cT_R(u^-;(x_0,t_0))
\end{split}\end{equation}
Taking $\dt=\f{1}{2c}$ in (7.38) and applying Young's inequality with $\vep=(b-a)^{(n+2s)(\f{1}{2}-\f{1}{p})}$ yield that
\begin{equation}\begin{split}
\sup_{Q^0_{ar}}u&\le\f{1}{2}\sup_{Q^0_{br}}u+\f{c}{(b-a)^{\f{n+2s}{p}}}\biggl(\f{1}{|\cQ^0_{2r}|}\iint_{Q^0_{2r}}u^p\,dx\,dt\biggr)^{\f{1}{p}} \\
&\qquad\qquad\qquad+c\biggl(\f{r}{R}\biggr)^{2s}\cT_R(u^-;(x_0,t_0)).
\end{split}\end{equation}
Employing Lemma 2.3 in (7.39) leads us to obtain that
\begin{equation}\begin{split}
\sup_{Q^0_{\vr r}}u&\le c\biggl[\f{1}{(2-\vr)^{\f{n+2s}{p}}}\biggl(\f{1}{|\cQ^0_{2r}|}\iint_{Q^0_{2r}}u^p\,dx\,dt\biggr)^{\f{1}{p}} +\biggl(\f{r}{R}\biggr)^{2s}\cT_R(u^-;(x_0,t_0))\biggr]
\end{split}\end{equation}
for any $p\in(0,2]$ and any $\vr\in[1/2,2)$.
We note that
\begin{equation}(A+B)^{\f{1}{p}}\le 2^{\f{1}{p}}(A^{\f{1}{p}}+B^{\f{1}{p}})
\end{equation} for any $A,B\ge 0$ and
\begin{equation}\f{2}{|\cQ^0_{2r}|}\le\f{1}{2|(Q^0_{2r})^+|}\,\,\Leftrightarrow\,\,\f{2}{2-\sm}\le\f{1}{2\sm}\,\,\Leftrightarrow\,\,0<\sm\le\f{2}{5}.
\end{equation}  Taking $\vr=1$ in (7.40) it follows from (7.41), (7.42) and Theorem 7.4 that
\begin{equation}\begin{split}
\sup_{Q^0_r}u&\le c\biggl(\f{1}{|\cQ^0_{2r}|}\iint_{Q^0_{2r}}u^p\,dx\,dt\biggr)^{\f{1}{p}} +c\biggl(\f{r}{R}\biggr)^{2s}\cT_R(u^-;(x_0,t_0)) \\
&\le c\biggl(\f{2}{|\cQ^0_{2r}|}\iint_{(Q^0_{2r})^+}u^p\,dx\,dt\biggr)^{\f{1}{p}} \\
&\quad+c\biggl(\f{2|Q^0_{2r}\s (Q^0_{2r})^+|}{|\cQ^0_{2r}|}\biggr)^{\f{1}{p}}\sup_{Q^0_{2r}\s (Q^0_{2r})^+}u+c\biggl(\f{r}{R}\biggr)^{2s}\cT_R(u^-;(x_0,t_0)) \\
&\le c\biggl(\inf_{(Q^0_r)^+}u+\f{4}{3}\biggl(\f{r}{R}\biggr)^{2s}\cT_R(u^-;(x_0,t_0))\biggr)\\
&\qquad\qquad\qquad\qquad\qquad+\f{1}{2}\,\sup_{(Q^0_r)^-}u+c\biggl(\f{r}{R}\biggr)^{2s}\cT_R(u^-;(x_0,t_0))
\end{split}\end{equation}
for a sufficiently small $p\in(0,2]$ (by (7.25) and (7.34)) satisfying
$$c\biggl(\f{2|Q^0_{2r}\s (Q^0_{2r})^+|}{|\cQ^0_{2r}|}\biggr)^{\f{1}{p}}\sup_{Q^0_{2r}\s (Q^0_{2r})^+}u=c\biggl(\f{2-2\sm}{2-\sm}\biggr)^{\f{1}{p}}\sup_{Q^0_{2r}\s (Q^0_{2r})^+}u
<\f{1}{2}\,\sup_{(Q^0_r)^-}u.$$
This implies that
\begin{equation}\begin{split}
\sup_{(Q^0_r)^-}u\le c\inf_{(Q^0_r)^+} u +c\biggl(\f{r}{R}\biggr)^{2s}\cT_R(u^-;(x_0,t_0)).
\end{split}\end{equation}
Hence we complete the proof. \qed

\section{Appendix: Existence and uniqueness of weak solutions}
The main goal of this appendix is to give the existence and uniqueness of weak solutions of the nonlocal parabolic equation ${\bf NE}_{\Om_I}(f,0,h)$ where $\ff\in L^2(I;L^2(\Om))$ and $h\in L^2(\Om)$ and even where  $\ff\in L^2(I;\rX^*_0(\Om))$ and $h\in L^2(\Om)$, and moreover is to give those of weak solutions of  the nonlocal parabolic equation ${\bf NE}_{\Om_I}(f,g,g)$  where  $\ff\in L^2(I;\rX^*_0(\Om))$ and $g\in H^s_T(\BR^n)$.

It is easy to check that the weak formulation of the following
eigenvalue problem
\begin{equation}\begin{cases}-\rL_K u=\ap\,u &\text{ in $\Om$ }\\\qquad\, u=0
&\text{ in $\BR^n\s\Om$, }
\end{cases}\end{equation} where $n>2s$, is given by
\begin{equation}\begin{cases}\la u,v\ra_{\rX_0(\Om)}=\ap\la u,v\ra_{L^2(\Om)},\,\forall v\in\rX_0(\Om),\\\qquad\,\,\,\,\, u\in\rX_0(\Om).
\end{cases}\end{equation}
Then it is well-known \cite{SV} that there exists a sequence
$\{\ap_i\}_{i\in\BN}$ of eigenvalues $\ap_i$ of (8.2) with
$0<\ap_1\le\ap_2\le\cdots\le\ap_i\le\ap_{i+1}\le\cdots$ and
$\lim_{i\to\iy}\ap_i=\iy$ such that the set $\{e_i\}_{i\in\BN}$ of
eigenfunctions $e_i$ corresponding to $\ap_i$ is an orthonormal
basis of $L^2(\Om)$ and an orthogonal basis of $\rX_0$. Moreover, it
turns out that $e_i\in\cQ_{i+1}$ and
\begin{equation}\ap_{i+1}=\| e_{i+1}\|^2_{\rX_0(\Om)}
\end{equation}
for any $i\in\BN$, where $\cQ_{i+1}=\{u\in\rX_0(\Om):\la
u,e_j\ra_{\rX_0(\Om)}=0,\,\forall j=1,2,\cdots,i\}$.

\,We construct a weak solution of the nonlocal
parabolic boundary value problem $(1.3)$ by using the eigenfunctions
of the nonlocal eigenvalue problem mentioned in (8.2), which is
called {\it Galerkin's approximation}.

\begin{thm}  If $\ff\in L^2(I;L^2(\Om))$ and $h\in L^2(\Om)$, then there exists
a weak solution $\fu\in H^1(I;\rX_0(\Om))$ of the nonlocal parabolic boundary value problem ${\bf NE}_{\Om_I}(f,0,h)$. Moreover,
$\fu\in C(I;\rX_0(\Om))$ after being modified on a set of measure zero.
\end{thm}

Let $\{e_i\}_{i\in\BN}$ be the set of eigenfunctions $e_i$
corresponding to eigenvalues $\ap_i$ of (8.2) that is an
orthonormal basis of $L^2(\Om)$ and an orthogonal basis of $\rX_0(\Om)$.
For $k\in\BN$, we consider a function $\fu_k:I\to\rX_0(\Om)$ of the form
\begin{equation}\fu_k(t)=\sum_{i=1}^k c^i_k(t)\,e_i,\,\,t\in I.
\end{equation}
Our next step is to show the existence of the functions
$\{c^i_k(t)\}_{i=1}^k$ for which
\begin{equation}c^i_k(-T)=\la h,e_i\ra_{L^2(\Om)}
\end{equation}
and
\begin{equation}\la\fu_k(t),e_i\ra_{\rX_0(\Om)}+\la\fu'_k(t),e_i\ra_{L^2(\Om)}
+\la\ff(t),e_i\ra_{L^2(\Om)}=0,
\,t\in I,
\end{equation} for any $i=1,2,\cdots,k$.

\begin{lemma}For each $k\in\BN$, there exists a unique function
$\fu_k$ of the form $(8.4)$ so that $(8.5)$ and $(8.6)$ hold.
\end{lemma}

\pf By applying (8.2) and (8.4), we reduce the weak formulation
(8.5) and (8.6) of (1.3) to the ordinary differential equations
\begin{equation}\begin{cases}\nu_i c^i_k(t)+{c_k^i}'(t)+\la\ff(t),e_i\ra_{L^2(\Om)}=0 \\
c^i_k(-T)=\la h,e_i\ra_{L^2(\Om)}
\end{cases}\end{equation} for any $i=1,2,\cdots,k$.
From standard O.D.E. theory, the initial value problem has a unique
solution $(c^1_k(t),c^2_k(t),\cdots,c^k_k(t))$ which satisfies (8.7)
for $\aee$ $t\in I$ and is absolutely continuous on $I$. Hence the
functions $\fu_k$ defined by (8.4) solves (8.5) and (8.6). \qed

Next we want to obtain a subsequence of the solutions $\fu_k$ of
(8.5) and (8.6) which converges to a weak solution of (1.3). To get
this, we need some uniform estimates which is called {\it energy
estimates}.

\begin{thm} If $\ff\in L^2(I;L^2(\Om))$ and $h\in L^2(\Om)$, then
the solutions $\fu_k$ obtained in Lemma 8.2 satisfy the following
energy estimates; that is, there exists a constant $c>0$ depending
only on $T,K$ and $\Om$ such that
\begin{equation*}\begin{split}\|\fu_k\|_{L^{\iy}(I;L^2(\Om))}&+\|\fu_k\|_{H^1(I;\rX_0(\Om))}
\le
c\bigl(\,\|\ff\|_{L^2(I;L^2(\Om))}+\|h\|_{L^2(\Om)}\bigr)\,\,\text{for all $k\in\BN$.}
\end{split}\end{equation*}
\end{thm}

\pf From (8.6), we can easily derive the equality
\begin{equation}\la\fu_k(t),\fu_k(t)\ra_{\rX_0(\Om)}+\la\fu'_k(t),\fu_k(t)\ra_{L^2(\Om)}+\la\ff(t),\fu_k(t)\ra_{L^2(\Om)}=0
\,\,\,\aee\,\,t\in I.
\end{equation}
Thus this yields the inequality
\begin{equation}\f{d}{dt}\,\|\fu_k(t)\|^2_{L^2(\Om)}+2\,\|\fu_k(t)\|^2_{\rX_0(\Om)}
\le\|\ff(t)\|^2_{L^2(\Om)}+\|\fu_k(t)\|^2_{L^2(\Om)}
\,\,\text{ for $\aee$ $t\in I$.}
\end{equation}
So it follows from the fractional Sobolev inequality and Gronwall's inequality that
\begin{equation}\begin{split}\|\fu_k(t)\|^2_{L^2(\Om)}&\le
e^{-c(t+T)}\biggl(\,\|h\|^2_{L^2(\Om)}+\int_{-T}^t\|\ff(\tau)\|^2_{L^2(\Om)}\,d\tau\biggr)\\
&\le\bigl(\,\|h\|^2_{L^2(\Om)}+\|\ff(t)\|^2_{L^2(I;L^2(\Om))}\bigr)
\end{split}\end{equation}
for $\aee$ $t\in I$; that is,
$$\|\fu_k\|^2_{L^{\iy}(I;L^2(\Om))}\le\bigl(\,\|h\|^2_{L^2(\Om)}+\|\ff(t)\|^2_{L^2(I;L^2(\Om))}\bigr).$$
By (8.9) and (8.10), we also obtain that
\begin{equation}\|\fu_k\|^2_{L^2(I;\rX_0(\Om))}=\int_{-T}^0\|\fu_k(\tau)\|^2_{\rX_0(\Om)}\,d\tau
\le c\bigl(\,\|h\|^2_{L^2(\Om)}+\|\ff(t)\|^2_{L^2(I;L^2(\Om))}\bigr)
\end{equation}

We see that $\rX_0(\Om)=\cP_k\oplus\cQ_{k+1}$ where
$\cP_k=\spa\{e_1,e_2,\cdots,e_k\}$ and $\cQ_{k+1}$ is the space given in
(8.3). Fix any $v\in\rX_0(\Om)$ with $\|v\|_{\rX_0(\Om)}\le 1$. Then we write
$v=v_1+v_2$ for $v_1\in\cP_k$ and $v_2\in\cQ_{k+1}$. Since
$\{e_i\}_{i\in\BN}$ are orthogonal in $\rX_0(\Om)$, we have that
$\|v_1\|_{\rX_0(\Om)}\le\|v\|_{\rX_0(\Om)}\le 1$. As in (8.6), we deduce that
\begin{equation}\la\fu_k(t),v_1\ra_{\rX_0(\Om)}+\la\fu'_k(t),v_1\ra_{L^2(\Om)}
+\la\ff(t),v_1\ra_{L^2(\Om)}=0
\end{equation}
for $\aee$ $t\in I$. Thus (8.4) and (8.12) imply that
$$\la\fu'_k(t),v\ra_{L^2(\Om)}=\la\fu'_k(t),v_1\ra_{L^2(\Om)}
=-\la\fu_k(t),v_1\ra_{\rX_0(\Om)}-\la\ff(t),v_1\ra_{L^2(\Om)}.$$
So it follows from Schwarz inequality on $\rX_0(\Om)$ that
$$|\la\fu'_k(t),v\ra_{L^2(\Om)}|\le
c\bigl(\,\|\fu_k(t)\|_{\rX_0(\Om)}+\|\ff(t)\|_{L^2(\Om)}\bigr).$$ This
gives that $$\|\fu'_k(t)\|_{\rX_0^*(\Om)}\le
c\bigl(\,\|\fu_k(t)\|_{\rX_0(\Om)}+\|\ff(t)\|_{L^2(\Om)}\bigr).$$ Therefore
by (8.11) we conclude that
\begin{equation}\|\fu'_k\|^2_{L^2(I;\rX_0^*)}\le
c\bigl(\,\|h\|^2_{L^2(D)}+\|\ff(t)\|^2_{L^2(D)}\bigr).
\end{equation}
Hence we are done.\qed

\,\,{\bf [Proof of Theorem 8.1.]} By Theorem 8.3 and Alaoglu's Theorem,
we see that there exist a subsequence
$\{\fu_{k_j}\}_{j\in\BN}\subset\{\fu_k\}_{k\in\BN}$ and a function
$\fu\in L^2(I;\rX_0(\Om))$ with $\fu'\in L^2(I;\rX_0^*(\Om))$ such that
\begin{equation}\begin{cases}
\fu_{k_j}\to\fu\,\,\text{ weakly in $L^2(I;\rX_0(\Om))$ }\\
\fu'_{k_j}\to\fu'\,\,\text{ weakly in $L^2(I;\rX_0^*(\Om))$ }
\end{cases}\text{ as $j\to\iy$. } \end{equation}
Let $C^1_c(I;\rX_0(\Om))$ be the set of all $\fv_N(t)\in C^1(I;\rX_0(\Om))$ of
the form $$\fv_N(t)=\sum_{i=1}^N c_i(t)e_i$$ where
$\{c_i\}_{i=1}^N\subset C^1_c(I)$ for some $N\in\BN$. Then we easily
see that $C^1_c(I;\rX_0(\Om))$ is dense in $L^2(I;\rX_0(\Om))$. Take any
$\fv_N(t)\in C^1(I;\rX_0(\Om))$ and choose $k\ge N$. Then it follows from
(8.6) and integrating with respect to $t$ that
\begin{equation}\int_{-T}^0\bigl[\la\fu_{k_j}(t),\fv_N(t)\ra_{\rX_0(\Om)}+\la\fu'_{k_j}(t),\fv_N(t)\ra_{L^2(\Om)}
+\la\ff(t),\fv_N(t)\ra_{L^2(\Om)}\bigr]\,dt=0.
\end{equation}
Passing (8.15) to weak limits (8.14) and using density of
$C^1_c(I;\rX_0(\Om))$ in $L^2(I;\rX_0(\Om))$, we obtain that
\begin{equation}\int_{-T}^0\bigl[\la\fu(t),\fv(t)\ra_{\rX_0(\Om)}+\la\fu'(t),\fv(t)\ra_{L^2(\Om)}
+\la\ff(t),\fv(t)\ra_{L^2(\Om)}\bigr]\,dt=0.
\end{equation}
for any $\fv(t)\in L^2(I;\rX_0(\Om))$. If we choose any elements of the
form $\fv(t)=\vp(t)v$ where $\vp\in C^1_c(I)$ and $v\in\rX_0(\Om)$, then
(8.16) becomes
\begin{equation}\int_{-T}^0\bigl[\la\fu(t),v\ra_{\rX_0(\Om)}+\la\fu'(t),v\ra_{L^2(\Om)}
+\la\ff(t),v\ra_{L^2(\Om)}\bigr]\vp(t)\,dt=0
\end{equation}
for any $\vp\in C^1_c(I)$. This implies that
$$\la\fu(t),v\ra_{\rX_0(\Om)}+\la\fu'(t),v\ra_{L^2(\Om)}
+\la\ff(t),v\ra_{L^2(\Om)}=0\,\,\text{ for $\aee$ $t\in I$ }$$ for any
$v\in\rX_0(\Om)$. Moreover, by the remark just above Definition 2.1, we see that $\fu\in C(I;L^2(\Om))$
after being modified on a set of measure zero.

In order to show that $\fu$ is a weak solution of (1.3), we finally
have only to prove that $\fu(-T)=h$. By (8.16) and integration by
parts, we have that
\begin{equation}\begin{split}&\int_{-T}^0\bigl[\la\fu(t),\fv(t)\ra_{\rX_0(\Om)}-\la\fu(t),\fv'(t)\ra_{L^2(\Om)}
+\la\ff(t),\fv(t)\ra_{L^2(\Om)}\bigr]\,dt\\&\qquad\qquad\qquad\qquad\qquad\qquad=-\la\fu(-T),\fv(-T)\ra_{L^2(\Om)}.
\end{split}\end{equation}
for any $\fv(t)\in L^2(I;\rX_0(\Om))$ with $\fv(0)=0$. Similarly it
follows from (8.15) and integration by parts that
\begin{equation}\begin{split}&\int_{-T}^0\bigl[\la\fu_{k_j}(t),\fv(t)\ra_{\rX_0(\Om)}-\la\fu_{k_j}(t),\fv'(t)\ra_{L^2(\Om)}
+\la\ff(t),\fv(t)\ra_{L^2(\Om)}\bigr]\,dt\\&\qquad\qquad\qquad\qquad\qquad\qquad=-\la\fu_{k_j}(-T),\fv(-T)\ra_{L^2(\Om)}.
\end{split}\end{equation}
for any $\fv(t)\in L^2(I;\rX_0(\Om))$ with $\fv(0)=0$. By (8.4) and (8.5), we get that
\begin{equation}\fu_{k_j}(-T)=\sum_{i=1}^{k_j}\la
h,e_i\ra_{L^2(\Om)}e_i. \end{equation} Since the orthonormal basis
$\{e_i\}_{i\in\BN}$ of $L^2(\Om)$ is complete, we see that
\begin{equation}\lim_{j\to\iy}\|\fu_{k_j}(-T)-h\|_{L^2(\Om)}=0.\end{equation}
Passing (8.19) to weak limits (8.14) and comparing it with (8.18)
and (8.21), we conclude that $\fu(-T)=h$. \qed

\,\,We shall obtain an improved regularity result for a weak solution of
(1.3) which is better than that of Theorem 8.1.

\begin{thm} Let $\ff\in L^2(I;L^2(\Om))$ and $h\in\rX_0(\Om)$. If $\,\fu\in
H^1(I;\rX_0(\Om))$ is a weak solution of
the nonlocal parabolic boundary value problem ${\bf NE}_{\Om_I}(f,0,h)$, then we have
that $\fu\in L^2(I;H^s(\Om))\cap
L^{\iy}(I;H^s(\Om))$ and $\fu'\in L^2(I;L^2(\Om))$, and
moreover there exists a constant $c_1>0$ depending only on $T,K$ and
$\Om$ such that
\begin{equation*}\begin{split}&\|\fu\|_{L^{\iy}(I;H^s(\Om))}+\|\fu\|_{L^2(I;H^s(\Om))}
+\|\fu'\|_{L^2(I;L^2(\Om))} \\
&\qquad\qquad\qquad\qquad\qquad\le
c_1\bigl(\,\|\ff\|_{L^2(I;L^2(\Om))}+\|h\|_{\rX_0(\Om)}\bigr).
\end{split}\end{equation*}
\end{thm}

\pf From (8.6), we have that for any $t\in I$,
\begin{equation*}\int_{-T}^t\biggl(\f{d}{dt}\,\|\fu_{k_j}(\tau)\|^2_{\rX_0(\Om)}+2\,\|\fu'_{k_j}(\tau)\|^2_{L^2(\Om)}
+2\,\la\ff(\tau),\fu'_{k_j}(\tau)\ra_{L^2(\Om)}\biggr)\,d\tau=0.
\end{equation*} From (8.20) and Bessel's inequality on the Hilbert space $\rX_0(\Om)$, we
can obtain that $$\|\fu_{k_j}(-T)\|^2_{\rX_0(\Om)}\le\|h\|^2_{\rX_0(\Om)}$$
for all $i\in\BN$. Applying Schwarz inequality on $L^2(\Om)$ and
Cauchy's inequality, this gives that
\begin{equation}\|\fu_{k_j}(t)\|^2_{\rX_0(\Om)}+\f{3}{2}\int_{-T}^t\|\fu'_{k_j}(\tau)\|^2_{L^2(\Om)}\,d\tau
\le\|h\|^2_{\rX_0(\Om)}+2\,\|\ff\|_{L^2(I;L^2(\Om))}
\end{equation}
for any $t\in I$. From (8.22), we easily get that
\begin{equation}\begin{split}&\|\fu_{k_j}\|^2_{L^{\iy}(I;\rX_0(\Om))}+\f{3}{2}\,\|\fu'_{k_j}\|^2_{L^2(I;L^2(\Om))}\\
&\qquad\qquad\qquad\qquad\qquad\le\|h\|^2_{\rX_0(\Om)}+2\,\|\ff\|_{L^2(I;L^2(\Om))}.
\end{split}\end{equation}
Since $C_c(I;\rX_0(\Om))$ is dense in $L^p(I;\rX_0(\Om))$ for all
$p\in[1,\iy)$, we note that
\begin{equation}\begin{split}&\|\fu_{k_j}\|_{L^{\iy}(I;\rX_0(\Om))}\\
&\,\,=\sup\biggl\{\int_{-T}^0\la\fu_{k_j}(t),\fh(t)\ra_{\rX_0(\Om)}\,dt:\|\fh\|_{L^1(I;\rX_0(\Om))}\le
1, \fh\in C_c(I;\rX_0(\Om))\biggr\} \end{split}\end{equation} by the
duality of  $L^1(I;\rX_0(\Om))$. So it follows from (8.14), (8.23) and
(8.24) that
\begin{equation}\int_{-T}^0\la\fu_{k_j}(t),\fh(t)\ra_{\rX_0(\Om)}\,dt\le\|h\|^2_{\rX_0(\Om)}+2\,\|\ff\|_{L^2(I;L^2(\Om))}
\end{equation} and thus, we have that
\begin{equation}\int_{-T}^0\la\fu(t),\fh(t)\ra_{\rX_0(\Om)}\,dt\le\|h\|^2_{\rX_0(\Om)}+2\,\|\ff\|_{L^2(I;L^2(\Om))}
\end{equation} for any $\fh\in C_c(I;\rX_0(\Om))$ with $\|\fh\|_{L^1(I;\rX_0(\Om))}\le
1$. Combining this with (2.6) and (2.7) imply that
\begin{equation}\|\fu\|^2_{L^{\iy}(I;H^s(\Om))}\le\|\fu\|^2_{L^{\iy}(I;\rX_0(\Om))}
\le\|h\|^2_{\rX_0(\Om)}+2\,\|\ff\|_{L^2(I;L^2(\Om))}.
\end{equation}

On the other hand, integrating (8.22) with respect to $t\in[-T,0)$,
we obtain that
\begin{equation}\begin{split}&\|\fu_{k_j}\|^2_{L^2(I;\rX_0(\Om))}+\f{3}{2}\,\|\fu'_{k_j}\|^2_{L^2(I;L^2(\Om))}\\
&\qquad\qquad\qquad\qquad\qquad\le\|h\|^2_{\rX_0(\Om)}+2\,\|\ff\|_{L^2(I;L^2(\Om))}.
\end{split}\end{equation}
By Riesz representation theorem, we know that
\begin{equation}\begin{split}&\|\fu'_{k_j}\|^2_{L^2(I;L^2(\Om))}\\
&=\sup\biggl\{\int_{-T}^0\la\fu'_{k_j}(t),\fy(t)\ra_{L^2(\Om)}\,dt:\|\fy\|_{L^2(I;L^2(\Om))}\le
1, \fy\in L^2(I;L^2(\Om))\biggr\}. \end{split}\end{equation} By (8.28)
and (8.29), we see that
\begin{equation}\int_{-T}^0\la\fu'_{k_j}(t),\fy(t)\ra_{L^2(\Om)}\,dt\le
c\bigl(\,\|h\|^2_{\rX_0(\Om)}+\,\|\ff\|_{L^2(I;L^2(\Om))}\bigr)
\end{equation} for all $\fy\in L^2(I;L^2(\Om))$ with $\|\fy\|_{L^2(I;L^2(\Om))}\le
1$. In particular, (8.30) holds for all $\fy\in L^2(I;\rX_0(\Om)$ with
$\|\fy\|_{L^2(I;L^2(\Om))}\le 1$. Thus by (8.14), we obtain that
\begin{equation}\begin{split}\lim_{j\to\iy}\int_{-T}^0\la\fu'_{k_j}(t),\fy(t)\ra_{L^2(\Om)}\,dt
&=\int_{-T}^0\la\fu'(t),\fy(t)\ra_{L^2(\Om)}\,dt\\
&\le c\bigl(\,\|h\|^2_{\rX_0(\Om)}+\,\|\ff\|_{L^2(I;L^2(\Om))}\bigr)
\end{split}\end{equation} for any $\fy\in L^2(I;\rX_0(\Om))$ with
$\|\fy\|_{L^2(I;L^2(\Om))}\le 1$. Also by Alaoglu's theorem and
(8.28), there are a subsequence
$\{\fu'_{k_{\el}}\}_{\el\in\BN}\subset\{\fu'_{k_j}\}_{j\in\BN}$ and
a function $\fv\in L^2(I;L^2(\Om))$ such that
\begin{equation*}\fu'_{k_{\el}}\to\fv\,\,\text{ weakly in
$L^2(I;L^2(\Om))$ }\,\text{ as $\el\to\iy$.}
\end{equation*} This and (8.30) yields that
\begin{equation}\begin{split}\lim_{\el\to\iy}\int_{-T}^0\la\fu'_{k_{\el}}(t),\fy(t)\ra_{L^2(\Om)}\,dt
&=\int_{-T}^0\la\fv(t),\fy(t)\ra_{L^2(\Om)}\,dt\\
&\le c\bigl(\,\|h\|^2_{\rX_0(\Om)}+\,\|\ff\|_{L^2(I;L^2(\Om))}\bigr)
\end{split}\end{equation} for all $\fy\in L^2(I;L^2(\Om))$ with
$\|\fy\|_{L^2(I;L^2(\Om))}\le 1$. Hence it follows from (8.31) and
(8.32) that
$$\int_{-T}^0\la\fu'(t),\fy(t)\ra_{L^2(\Om)}\,dt=\int_{-T}^0\la\fv(t),\fy(t)\ra_{L^2(\Om)}\,dt$$
for any $\fy\in L^2(I;\rX_0(\Om))$ with $\|\fy\|_{L^2(I;L^2(\Om))}\le 1$.
This implies $\fu'=\fv$ $\aee\,$. Thus by (8.32) we get that
$\fu'\in L^2(I;L^2(\Om))$ and
$$\int_{-T}^0\la\fv(t),\fy(t)\ra_{L^2(\Om)}\,dt\le
c\bigl(\,\|h\|^2_{\rX_0(\Om)}+\|\ff\|_{L^2(I;L^2(\Om))}\bigr)$$ for all
$\fy\in L^2(I;L^2(\Om))$ with $\|\fy\|_{L^2(I;L^2(\Om))}\le 1$. Hence we
conclude that
$$\|\fu'\|_{L^2(I;L^2(\Om))}\le
c\bigl(\,\|h\|^2_{\rX_0(\Om)}+\,\|\ff\|_{L^2(I;L^2(\Om))}\bigr).$$

Finally, we can derive the estimate
$$\|\fu\|^2_{L^2(I;H^s(\Om))}\le\|\fu\|^2_{L^2(I;\rX_0(\Om))}
\le c\bigl(\,\|h\|^2_{\rX_0(\Om)}+\,\|\ff\|_{L^2(I;L^2(\Om))}\bigr)$$ from
(2.6), (2.7), (8.14) and (8.28). Hence we complete the proof. \qed

\,\,Even when $\ff\in L^2(I;\rX^*_0(\Om))$ and $h\in L^2(\Om)$, using duality we can obtain the following estimates in the similar way as in Theorem 8.3.

\begin{thm} If $\,\ff\in L^2(I;\rX^*_0(\Om))$ and $h\in L^2(\Om)$, then
the solutions $\fu_k$ obtained in Lemma 8.2 satisfy the following
energy estimates; that is, there exists a constant $c>0$ depending
only on $T,K$ and $\Om$ such that
\begin{equation*}\begin{split}\|\fu_k\|_{L^{\iy}(I;L^2(\Om))}&+\|\fu_k\|_{H^1(I;\rX_0(\Om))}
\le
c\bigl(\,\|\ff\|_{L^2(I;\rX^*_0(\Om))}+\|h\|_{L^2(\Om)}\bigr)\,\,\text{for all $k\in\BN$.}
\end{split}\end{equation*}
\end{thm}

\,\,By applying Theorem 8.5, we can prove the following theorem in the similar way as in Theorem 8.1.

\begin{thm}  If $\,\ff\in L^2(I;\rX^*_0(\Om))$ and $h\in L^2(\Om)$, then there exists
a weak solution $\fu\in H^1(I;\rX_0(\Om))$ of the nonlocal parabolic boundary value problem ${\bf NE}_{\Om_I}(f,0,h)$. Moreover,
$\fu\in C(I;\rX_0(\Om))$ after being modified on a set of measure zero.
\end{thm}

\begin{thm}  If $\,\ff\in L^2(I;\rX^*_0(\Om))$ and $g\in H^s_T(\BR^n)$, then there is
a weak solution $\fu\in H^1(I;\rX_g(\Om))$ of the nonlocal parabolic boundary value problem ${\bf NE}_{\Om_I}(f,g,g)$. Moreover,
$\fu-\fg\in C(I;\rX_0(\Om))$ after being modified on a set of measure zero.
\end{thm}

\pf By Theorem 8.6, there is a weak solution $\fv\in H^1(I;\rX_0(\Om))$ of the nonlocal parabolic boundary value problem ${\bf NE}_{\Om_I}(f-\rL_K g-\pa_t g,0,0)$, because it is easy to check that $f-\rL_K g-\pa_t g\in L^2(I;\rX^*_0(\Om))$. Then we easily see that $\fu=\fv+\fg$ satisfies the equation ${\bf NE}_{\Om_I}(f,g,g)$ and $\fu\in H^1(I;\rX_g(\Om))$. Hence we are done. \qed

\,\,In the following theorem, we give the uniqueness of the weak solutions of the nonlocal parabolic equation given in (1.3). Its proof is quite simple and follows from a direct application of Gronwall's inequality.

\begin{thm} A weak solution $\fu\in H^1(I;\rX_g(\Om))$ of the nonlocal parabolic boundary value
problem ${\bf NE}_{\Om_I}(f,g,h)$ is unique, if it exists.
\end{thm}

\pf We have only to check that the only weak solution of the equation ${\bf NE}_{\Om_I}(f,g,h)$ with
$\ff=\bold 0$ and $g=h=0$ must be $\fu\equiv\bold 0\,\,\aee\,.$ By
(2.12), we have that
$$\f{d}{dt}\,\|\fu(t)\|^2_{L^2(\Om)}\le\f{d}{dt}\,\|\fu(t)\|^2_{L^2(\Om)}+2\|\fu(t)^2\|_{\rX_0(\Om)}=0$$ for $\aee$
$t\in I$. Thus it follows from Gronwall's inequality that
$$\|\fu(t)\|^2_{L^2(\Om)}\le\|g\|_{L^2(\Om)}=0.$$
Hence we conclude that $\fu\equiv\bold 0$ $\aee\,$. \qed

\,\,\,\noindent{\bf Acknowledgement.} Yong-Cheol Kim was supported by the National Research Foundation of Korea(NRF) grant funded by the Korea government (MSIP-2017R1A2\\B1005433).

\end{document}